\documentclass[12pt,a4paper,reqno]{amsart}
\usepackage{amsmath}
\usepackage{physics}
\usepackage{amssymb}
\numberwithin{equation}{section}
\usepackage{graphicx}
\usepackage{amssymb}
\usepackage{amsthm}
\usepackage{enumitem}
\usepackage{tikz}
\usetikzlibrary{shapes}
\usepackage{mathrsfs}
\usepackage{hyperref}
\usetikzlibrary{arrows}
\usepackage[utf8]{inputenc}
\usetikzlibrary{calc}
\numberwithin{equation}{section}
\usepackage{tikz-cd}
\usepackage[utf8]{inputenc}
\usepackage[T1]{fontenc}

\usepackage{tikz-cd}
\usepackage{pgfplots}
\pgfplotsset{compat=1.10}
\usepgfplotslibrary{fillbetween}
\usepackage{txfonts}
\usepackage{lipsum} 
\let\oldpart\part
\renewcommand\part{\newpage\oldpart}

\usepackage[%
shortalphabetic, 
abbrev, 
backrefs, 
]{amsrefs}
\hypersetup{%
	linktoc=all, 
	colorlinks=true,
	citecolor=magenta,
}

\NewDocumentCommand\Crefnameitem { m m m O{\textup} O{(\roman*)}} {%
	\Crefname{#1enumi}{#2}{#3} 
	\AtBeginEnvironment{#1}{%
		\crefalias{enumi}{#1enumi}%
		\setlist[enumerate,1]{
			label={#4{#5}.},
			ref={#5}
		}%
	}  
}

\newtheorem{thm}{Theorem}[section]
\newtheorem*{theorem*}{Theorem}

\newtheorem{prop}[thm]{Proposition}
\newtheorem{lm}[thm]{Lemma}

\newtheorem{coro}[thm]{Corollary}

\newcommand\numberthis{\addtocounter{equation}{1}\tag{\theequation}}

\newtheorem{innercustomgeneric}{\customgenericname}
\providecommand{\customgenericname}{}
\newcommand{\newcustomtheorem}[2]{%
	\newenvironment{#1}[1]
	{%
		\renewcommand\customgenericname{#2}%
		\renewcommand\theinnercustomgeneric{##1}%
		\innercustomgeneric
	}
	{\endinnercustomgeneric}
}
\newcustomtheorem{customthm}{Theorem}

\theoremstyle{definition}

\newtheorem{df}[thm]{Definition}
\newtheorem{remark}[thm]{Remark}

\allowdisplaybreaks

\newcommand{\N}{\mathbb{N}}
\newcommand{\Z}{\mathbb{Z}}

\newcommand{\R}{\mathbb{R}}

\def \ra {\rightarrow}

\def \C {\mathbb{C}}

\def\la{\lambda}

\def \al{\alpha}

\def \om{\omega}
\def \ga {\gamma}

\def \b {\beta}
\def \op {\oplus}

\def \ssq{\subseteq}
\def \vac {\mathbf{1}}

\def \h {\mathfrak{h}}

\def \End {\mathrm{End}}

\def\Id{\mathrm{Id}}
\def \wt {\mathrm{wt}}
\def\ep{\epsilon}

\def \bs {\backslash}

\def\L{\mathcal{L}}

\def\ds{\dots}
\def\o{\otimes}

\def\spn{\mathrm{span}}
\def\cone{\mathrm{Cone}}
\def\Om{\Omega}

\def\tY{\tilde{Y}}

\def\scrU{\mathscr{U}}

\def\<{\langle}
\def\>{\rangle}

\def\rN{\mathrm{N}}
\def\L{\mathsf{L}}

\def\fA{\mathfrak{A}}

\def\rN{\mathrm{N}}
\def\L{\mathsf{L}}
\def\sR{\mathsf{R}}

\def\A{\mathsf{A}}

\begin {document}

\title{Rank-two parabolic-type VOAs and nilpotency of nil ideals}


\author{Jianqi Liu }
\address{Department of Mathematics, University of Pennsylvania, Philadelphia, PA, USA, 19104}
\email{jliu230@sas.upenn.edu}

\maketitle

\begin{abstract}
In this paper, we undertake a systematic study of the parabolic-type sub-vertex operator algebras (subVOAs) \(V_P\) of rank-two lattice VOAs \(V_L\), originally introduced by the first-named author. We first classify all possible types of such subVOAs by analyzing the corresponding submonoids \(P \subseteq L\). For each type of \(V_P\), we then classify its irreducible modules. Certain Zhu algebras \(A(V_P)\) provide new examples of rings with nil ideals that are not nilpotent. Finally, we show that the simple quotient \(V_H\) of any parabolic-type subVOA \(V_P\) is a \(C_1\)-cofinite irrational VOA satisfying the strongly unital property recently introduced by Damiolini--Gibney--Krashen.
\end{abstract}

\tableofcontents

\allowdisplaybreaks

\section{Introduction}
Parabolic-type subVOAs of lattice VOAs, introduced by the first-named author \cite{Liu25}, arise naturally in the quasi-triangular decomposition of lattice VOAs and exhibit many properties analogous to parabolic subalgebras of semisimple Lie algebras. These subVOAs also provide natural solutions to the operator-form classical Yang--Baxter equations on VOAs \cite{BGL,BGLW}. In \cite{Liu25(1)}, a standard parabolic-type subVOA $V_P$ of the rank-two lattice VOA $V_{A_2}$ was used as a motivating example to study degree-zero induction for the embedding $V_P \hookrightarrow V_{A_2}$. In this paper, we develop a complete structural and representation-theoretic theory for all parabolic-type subVOAs of rank-two lattice VOAs, generalizing the rank-one results \cite[Section 6]{Liu25} as well as the special rank-two type-A case in \cite{Liu25(1)}. Our methods naturally extend to higher-rank lattices. Furthermore, the simple quotients $V_H$ of these rank-two parabolic-type VOAs form a natural class of CFT-type $C_1$-cofinite irrational VOAs whose mode transition algebras $\mathfrak{A}_d$, recently introduced by Damiolini--Gibney--Krashen \cite{DGK23}, satisfy the strongly unital property.

The lattice vertex operator algebra $V_L$, introduced by Frenkel--Lepowsky--Meurman \cite{FLM}, is the first and most fundamental example of a VOA and plays a central role in the theory. For example, twisted modules and orbifold subVOAs of the Leech lattice VOA $V_\Lambda$ led to the original construction of the moonshine module $V^\natural$ \cite{FLM}, and $\Z_2$- and $\Z_3$-orbifolds of lattice VOAs were crucial in the proof of Schellekens' conjecture \cite{S93} on holomorphic VOAs of central charge $24$ \cite{EMS,Lam11,LS15,LS16}. Sub-lattice VOAs in affine VOAs also play key roles in the representation theory and rationality of parafermion VOAs \cite{DR17}. 

Although the representation theory of lattice VOAs is well understood \cite{D,D94,DL}, their structural theory continues to yield interesting results. In studying solutions to the classical Yang--Baxter equations on VOAs \cite{BGL,BGLW}, we observed that a lattice VOA $V_L = \bigoplus_{\alpha \in L} M_{\hat{\h}}(1,\alpha)$ is structurally analogous to a semisimple Lie algebra if we regard $L$ as an analog of the root system and the Heisenberg modules $M_{\hat{\h}}(1,\alpha)$ as an analog of the Cartan subalgebra. By taking submonoids $M \subset L$, one obtains many subVOAs $V_M = \bigoplus_{\alpha \in M} M_{\hat{\h}}(1,\alpha)$ of $V_L$. Some of these subVOAs are particularly notable. For instance, a subVOA $V_B \leq V_{A_2}$ of this kind provides an example of a CFT-type, $C_1$-cofinite subVOA whose Zhu algebra $A(V_B)$ is not noetherian \cite[Section 5.2]{Liu25}; we call such a VOA a Borel-type subVOA of $V_{A_2}$, a special type of rank-two parabolic-type VOA. Another standard rank-two parabolic-type subVOA $V_P \leq V_{A_2}$ served as a motivating example for the definition of degree-zero induction for VOA embeddings \cite{Liu25(1)}, and its Zhu algebra $A(V_P)$ was shown to be a nilpotent extension of a skew-polynomial ring. Motivated by these intriguing properties, we provide a systematic study of rank-two parabolic-type VOAs in this paper, extending the previous work of the first-named author in \cite[Section 6]{Liu25}.

To state our main results and outline their proofs, we first introduce some notation. 
Analogous to a Borel subalgebra (resp. Borel subgroup) in a semisimple Lie algebra (resp. a reductive algebraic group), a \emph{Borel-type submonoid} contains half of the ``root spaces.'' More precisely, it is a submonoid $B \leq L$ such that $B \cup (-B) = L$ and $B \cap (-B) = \{0\}$. One can show that there exists a unique hyperplane $P(\gamma)$ in the Euclidean space $E$ spanned by $L$ such that $B$ is contained in the positive half-space defined by $P(\gamma)$.  

We say that a submonoid $P \leq L$ is of \emph{parabolic-type} if it contains a Borel-type submonoid $B \leq L$. A parabolic-type subVOA of $V_L$ is then a subVOA 
\[
V_P = \bigoplus_{\alpha \in P} M_{\hat{\h}}(1, \alpha)
\] 
associated with a parabolic-type submonoid $P \leq L$. In particular, the parabolic-type subVOAs of a rank-two lattice VOA $V_L$ are completely determined by the corresponding parabolic-type submonoids.  
The following theorem provides the classification of parabolic-type submonoids in the rank-two case; see Theorem~\ref{thm:classificationofparamonoid} for details.  

\begin{customthm}{A}\label{main:A}
	Let $P$ be a parabolic-type submonoid of a rank-two even lattice $L$. Then it has one of the following two types: 
	\begin{enumerate}
		\item $P = B$ is a rank-two Borel-type submonoid. In this case, 
		\[
		P = (\Z_{\ge 0} \alpha) \cup (P^+(\gamma) \cap L),
		\] 
		where $P(\gamma)$ is the unique hyperplane such that $B \subset P^{\ge 0}(\gamma)$ and $P(\gamma) \cap L = \Z \alpha$ for some $\alpha \in L$ (possibly zero).
		\item $P = P^{\ge 0}(\gamma) \cap L$ for some hyperplane $P(\gamma) \subset E$ such that $P(\gamma) \cap L = \Z \alpha \neq \{0\}$.
	\end{enumerate}
	See Figures~\ref{fig2} and \ref{fig3} for an illustration in a standard rank-two lattice. We refer to $P$ as \emph{type-I} if $P = B$, and as \emph{type-II} if $P = P^{\ge 0}(\gamma) \cap L$ with $P(\gamma) \cap L \neq \{0\}$. 
\end{customthm}

The proof of Theorem~\ref{main:A} involves a technical argument regarding the existence of a basis of the rank-two lattice $L$ contained in a particular submonoid, see Proposition~\ref{prop:monoidhyperplane}. In analogy with semisimple Lie algebras, where a parabolic subalgebra $\mathfrak{p} \leq \mathfrak{g}$ contains a ``Cartan-part'' $\mathfrak{h}$. For example, in type $A$, $\mathfrak{h}$ is generated by block-diagonal matrices. The representation theory of $\mathfrak{p}$ is governed by $\mathfrak{h}$. Our second main theorem shows that a similar Cartan-part subVOA $V_H$ exists in any rank-two parabolic-type subVOA $V_P$, further justifying the terminology of ``parabolic-type subVOA''; see Proposition~\ref{prop:VHtensor} and Theorem~\ref{thm:strucVP}.

\begin{customthm}{B}\label{main:B}
	Let $V_P$ be a rank-two parabolic-type VOA. Then there exists a unique maximal proper ideal $V^+ \lhd V_P$ and a simple subVOA $V_H \leq V_P$ such that 
	\[
	V_P = V_H \oplus V^+ \quad \text{as vector spaces,} \quad \text{and} \quad V_P / V^+ \cong V_H \quad \text{as VOAs.}
	\] 
	Moreover, the simple quotient $V_H$ admits the following characterization: 
	\begin{enumerate}
		\item If $P = B$ is of type-I, then $V_H = M_{\hat{\h}}(1,0)$ is the rank-two Heisenberg VOA. 
		\item If $P$ is of type-II, with $P(\gamma) \cap L = \Z \alpha$, then 
		\[
		V_H = \bigoplus_{n \in \Z} M_{\hat{\h}}(1,n\alpha) \cong M_{\widehat{\C \beta}}(1,0) \otimes V_{\Z \alpha},
		\] 
		where $V_{\Z \alpha}$ is the rank-one lattice VOA associated to $\Z \alpha$, $\beta \in \h = \C \otimes_\Z L$ satisfies $\C \beta = (\C \alpha)^\perp$, and $M_{\widehat{\C \beta}}(1,0)$ is the rank-one Heisenberg VOA associated to $\C \beta$.
	\end{enumerate}
	We refer to $V_H$ as the \emph{Cartan-part} of $V_P$.
\end{customthm}

It is straightforward to see that $V_P$ admits a split decomposition $V_P = V_H \oplus V^+$ into a simple subVOA and an ideal. However, in the type-II case, the isomorphism 
\[
V_H \cong M_{\widehat{\C \beta}}(1,0) \otimes V_{\Z \alpha}
\] 
is not immediately obvious. This structure was strongly suggested by computations of the Zhu algebra for the standard rank-two type-$A$ parabolic-type subVOA $V_P$ in \cite{Liu25(1)}. To rigorously establish the isomorphism of VOAs, we carefully analyze the behavior of lattice vertex operators, as defined by Frenkel-Lepowsky-Meurman \cite{FLM}, under the decomposition 
\(\h = \C \otimes_\Z L = \C \alpha \oplus (\C \alpha)^\perp\), see Proposition~\ref{prop:VHtensor}.  
With the explicit structural description of $V_P$ in Theorem~\ref{main:B}, one can readily verify the $C_1$-cofiniteness of $V_H$ and, under mild conditions on the lattice $L$ \eqref{eq:conditionforranktwoC_1}, also of $V_P$ itself; see Theorem~\ref{thm:C1forVP}. In particular, $V_H$ is always $C_1$-cofinite, while $V_P$ is $C_1$-cofinite only under the specified condition.

Similar to the situation in classical Lie theory, the representation theory of $V_P = V_H \oplus V^+$ is governed by its Cartan part $V_H$. In particular, the maximal ideal $V^+$ acts trivially on any irreducible $V_P$-module $M$, and the irreducible $V_P$-modules are in one-to-one correspondence with irreducible $V_H$-modules. However, given an irreducible $V_P$-module $W$, it is not immediately clear that $V^+$ acts as zero on $ W$. 
To establish this, we employ the crucial tool of Zhu's algebra $A(V)$ \cite{Z}. We show that the image $A(V^+)$ forms a nil ideal in $A(V_P)$, which then acts as zero on any irreducible $A(V_P)$-module $W(0)$, see Proposition~\ref{prop:nilideal} and Lemma~\ref{lm:irrvanishing}. Conversely, given an irreducible $V_H$-module $(M, Y_M)$, one can extend it to an irreducible $V_P$-module $(M, \tilde{Y}_M)$ by defining $\tilde{Y}_M|_{V^+} = 0$. The following theorem provides a complete classification of irreducible $V_P$-modules, extending the rank-two type-$A$ case in \cite{Liu25(1)}. 

\begin{customthm}{C}\label{main:C}
	Let $V_P = V_H \oplus V^+$ be a rank-two parabolic-type VOA, and let $M = \bigoplus_{n=0}^{\infty} M(n)$ be an irreducible admissible $V_P$-module. Then $M$ is an irreducible $V_H$-module on which $V^+$ acts trivially. In particular, every irreducible admissible $V_P$-module is ordinary. Furthermore,
	\begin{enumerate}
		\item If $V_P$ is of type-I, then 
		\[
		\{ (M_{\hat{\h}}(1,\lambda), \tilde{Y}_M) : \lambda \in \h \}
		\] 
		form a complete list of irreducible $V_P$-modules up to isomorphism.
		\item If $V_P$ is of type-II, with $P(\gamma) \cap L = \Z \alpha$, $(\alpha|\alpha) = 2N$, and $\C \beta = (\C \alpha)^\perp$, then
		\[
		\Bigl\{ \bigl(L^{(\mu,i)} = M_{\widehat{\C \beta}}(1,\mu) \otimes V_{\Z \alpha + \frac{i}{2N} \alpha}, \tilde{Y}_M \bigr) : \mu \in \C \beta,\, i = 0,1,\dots,2N-1 \Bigr\}
		\] 
		form a complete list of irreducible $V_P$-modules up to isomorphism.
	\end{enumerate}
\end{customthm}


In a recent study of the smoothing property of conformal blocks associated with VOA modules, Damiolini-Gibney-Krashen \cite{DGK23} introduced a sequence of associative algebras, called the {\em mode transition algebras} $\fA_d$ for $d \in \N$. These algebras are defined as tensor products of certain subquotients of the universal enveloping algebra $\scrU(V)$ of a VOA $V$, and they satisfy the property that $\fA_d$ surjects onto the kernel of the canonical map $\A_d \to \A_{d-1}$ between higher-level Zhu's algebras \cite{DLM2} (see Definition~\ref{df:MTA}). It was shown in \cite[Remark 3.4.6]{DGK23} that if $V$ is rational, then $\fA_d$ possesses a {\em strong identity} for all $d \in \N$, which is equivalent to the smoothing property for VOA-conformal blocks \cite[Theorem 5.0.3]{DGK23}. However, rationality is not necessary for the existence of strong identity in $\fA_d$; for example, it was proved in \cite{DGK24} that the Heisenberg VOA satisfies the strongly unital condition for mode transition algebras.
Since rationality is preserved under tensor products of VOAs \cite{DMZ94} and the Cartan-part subVOA $V_H$ of a type-II parabolic-type VOA $V_P$ is a tensor product of two VOAs, $M_{\widehat{\C \beta}}(1,0)$ and $V_{\Z \alpha}$, each satisfying the strongly unital condition, it is natural to expect that $V_H$ also satisfies this property. Our final main theorem confirms this expectation.

\begin{customthm}{D}\label{main:D}
	Let $V_P = V_H \oplus V^+$ be a rank-two parabolic-type VOA. Then its Cartan-part subVOA $V_H$ satisfies the strongly unital condition for mode transition algebras.
\end{customthm}

This paper is organized as follows. In Section~\ref{sec:2}, we recall the definitions and basic properties of conic, Borel, and parabolic-type submonoids of an even lattice, and we prove the classification theorem for rank-two parabolic-type submonoids (Theorem~\ref{main:A}). In Section~\ref{sec:3}, we review the definition of lattice VOAs and establish the structural properties of rank-two parabolic-type subVOAs $V_P$ (Theorem~\ref{main:B}). Section~\ref{sec:4} is devoted to the classification of irreducible $V_P$-modules and the determination of fusion rules among them (Theorem~\ref{main:C}). Finally, in Section~\ref{sec:5}, we prove that the mode transition algebras $\fA_d$ of the Cartan-part subVOA $V_H$ of type-II parabolic-type VOAs are strongly unital for all $d \ge 0$ (Theorem~\ref{main:D}). Throughout, unless otherwise stated, all vector spaces are defined over $\C$.

\section{The rank-two parabolic-type monoids}\label{sec:2}

In this Section, we recall the notions of Borel and parabolic-type submonoids of a positive-definite even lattice $L$ introduced by the first named author in \cite{Liu25}. Then we classify all the rank-two parabolic-type monoids as a preparation for the study of rank-two parabolic-type VOAs.

\subsection{Basics of parabolic-type submonoids of an even lattice}

Let $\gamma$ be a nonzero vector in a Euclidean space $(E,(\cdot|\cdot))$. Let $P(\gamma)$ denote the hyperplane passing through the origin and perpendicular to $\gamma$. We define the positive and non-negative sides of $P(\gamma)$ by $P^+(\gamma)$ and $P^{\ge 0}(\gamma)$, respectively, and the opposite sides by $P^-(\gamma)$ and $P^{\le 0}(\gamma)$. More precisely,
\begin{align*}
	P(\gamma) &= \{ v \in E : (\gamma|v) = 0 \},\\
	P^+(\gamma) &= \{ v \in E : (\gamma|v) > 0 \}, \\
	P^{\ge 0}(\gamma) &= P(\gamma) \sqcup P^+(\gamma) = \{ v \in E : (\gamma|v) \ge 0 \}, \numberthis \label{eq:hyperplanes} \\
	P^-(\gamma) &= -P^+(\gamma) = \{ v \in E : (\gamma|v) < 0 \},\\
	P^{\le 0}(\gamma) &= -P^{\ge 0}(\gamma) = \{ v \in E : (\gamma|v) \le 0 \}.
\end{align*}

An \emph{even lattice} $L$ in $(E,(\cdot|\cdot))$ is a free abelian group of rank $r = \dim_\R E$ such that $(\alpha|\alpha) \in 2\mathbb{Z}$ for all $\alpha \in L$. A \emph{submonoid} $M \le L$ is an additive sub-abelian group of $L$ containing $0$.

\begin{df}\label{def:submonoidsofaPDElattice}\cite[Definition 2.1]{Liu25}
	Let $L$ be an even lattice in a Euclidean space $(E,(\cdot| \cdot))$.
	\begin{enumerate}
		\item A submonoid $C \le L$ is said to be of \textbf{conic-type} if there exists a basis $\{\alpha_1, \dots, \alpha_r\}$ of $L$ such that
		\[
		C = \Z_{\ge 0} \alpha_1 \oplus \dots \oplus \Z_{\ge 0} \alpha_r.
		\]
		\item A submonoid $B \le L$ is said to be of \textbf{Borel-type} if it satisfies:
		\begin{enumerate}
			\item $B \cup (-B) = L$,
			\item $B \cap (-B) = \{0\}$,
			\item There exists a hyperplane $P(\gamma) \subset E$ such that $B \subset P^{\ge 0}(\gamma)$.
		\end{enumerate}
		\item A proper submonoid $P \le L$ is said to be of \textbf{parabolic-type} if there exists a Borel-type submonoid $B \le L$ such that $B \subset P$.
	\end{enumerate}
\end{df}

Recall the following basic facts about Borel and parabolic-type monoids in  \cite[Propositions 2.2, 2.4]{Liu25}. 
\begin{prop}\label{prop:basicBorel}
	Let $L$ be an even lattice in the Euclidean space $(E,(\cdot| \cdot))$. Then 
	\begin{enumerate}
		\item Borel and parabolic-type submonoids of $L$ exist. They are invariant under the action of $\mathrm{Aut}(L)$.
		\item The hyperplane $P(\ga)$ in the definition of Borel-type submonoid is unique. In other words, given a Borel-type submonoid $B$, with $B\subset P^{\geq 0}(\ga)$, if there exists another $\ga'\in E\bs\{0\}$ such that $B\subset P^{\geq 0}(\ga')$, then $P(\ga)=P(\ga')$. 
	\end{enumerate}
\end{prop}


From now on, we fix a rank-two even lattice $L=\Z\al_1\op \Z\al_2$ in a two-dimensional Euclidean space $(E,(\cdot|\cdot))$. The following notion is useful for our later discussion.

\begin{df}\label{def:basiscone}
	Let $L$ be a rank-two even lattice, and  $\al_1,\al_2\in L\bs\{0\}$. Write
	\begin{equation}
		\cone(\al_1,\al_2):=\Z_{\geq 0}\al_1+\Z_{\geq 0}\al_2=\{m_1\al_1+m_2\al_2:m_1,m_2\in \Z_{\geq 0}\}. 
	\end{equation}
	We call $\cone(\al_1,\al_2)$ a {\bf basis cone} if $\{\al_1,\al_2\}$ is a basis of the lattice $L$. In this case, $C=\cone(\al_1,\al_2)$ is the conic-type submonid spanned by $\al_1$ and $\al_2$. 
	
\end{df}
The following properties about $\cone(\al_1,\al_2)$ are obvious, we omit the proof. 

\begin{lm}\label{lm:propertiesofcone}
	Let $L$ be a rank-two even lattice, and $\al_1,\al_2\in L\bs\{0\}$ be $\Z$-linearly independent. 
	\begin{enumerate}
		\item If $M\ssq L$ is a submonoid that contains $\al_1,\al_2$, then $\cone(\al_1,\al_2)\ssq M$.
		\item If $\cone(\al_1,\al_2)$ is a basis cone, then it contains all the lattice points of $L$ in the conic area $\R_{\geq 0}\al_1+\R_{\geq 0}\al_2$ of the Euclidean space $E$.
		\item  If $\{\al_1,\al_2\}$ a basis of $L$, then $L$ is a union of basis cones: 
		\begin{equation}\label{eq:coneforL}
			L=\cone(\al_1,\al_2)\cup \cone(\al_1,-\al_2)\cup \cone(-\al_1,\al_2)\cup \cone(-\al_1,-\al_2). 
		\end{equation}
	\end{enumerate}
\end{lm}

\begin{lm}\label{lm:basis}
	The submonoid $M=P^+(\ga)\cap L$ contains a basis $\{\al_1,\al_2\}$ of $L$.
\end{lm}
\begin{proof}
	Let $\{\b_1,\b_2\}$ be any basis of $L$. Note that $\b_1,\b_2$ cannot be both lying on $P(\ga)$, since otherwise the rank of $L$ would be one. If $\b_1,\b_2\in P^{\pm}(\ga)$, then the basis $\{\al_1=\pm\b_1,\al_2=\pm\b_2\}$ is contained in $P^+(\ga)\cap L$ of $L$, we are done. Otherwise, we may assume $\b_1\in P^+(\ga)$ and $\b_2\in P^{\leq 0}(\ga)$. Then the basis $\{\al_1=\b_1,\al_2=\b_1-\b_2\}$ is contained in $ P^+(\ga)\cap L$. 
\end{proof}

\begin{df}
	A nonzero lattice point $\al\in L\bs\{0\}$ is called {\bf primitive} if $\al$ is not a multiple of any other lattice points in $L$. i.e., if $\al=k\b$ for some $k\in \Z$ and $\b\in L$, then $k=\pm 1$. 
\end{df}

\subsection{Classification of the rank-two parabolic-type submonoids}
We need the following elementary result about the rank-two lattices for our classification theorem.
\begin{prop}\label{prop:monoidhyperplane}
	Let $P(\ga)$ be a hyperplane in $E$, $\al\in P^{-}(\ga)\cap L$ be a primitive lattice point, and $M$ be the submonoid in $L$ such that $\al\in M$ and $P^+(\ga)\cap L\subset M$. Then $M=L$.
\end{prop}
\begin{proof}
	By Lemma~\ref{lm:basis}, there exists a basis $\{\al_1,\al_2\}$ of $L$ that is contained in $ P^+(\ga)\cap L\subset M$. Let 
	\begin{equation}\label{eq:condition}
		(\al_1|\ga)=a_1>0,\quad (\al_2|\ga)=a_2>0,\quad \mathrm{and}\quad \al=m\al_1+n\al_2.
	\end{equation}
	Note that $(\al|\ga)=ma_1+na_2<0$. If $m<0$ and $n<0$, then 
	$$-\al_1=\al+(-m-1)\al_1+(-n)\al_2\in M,\quad \mathrm{and}\quad -\al_2=\al+(-m)\al_1+(-n-1)\al_2\in M.$$
	It follows from Lemma~\ref{lm:propertiesofcone} that $\cone(\pm \al_1,\pm\al_2)\ssq M$, and so $L=M$ by \eqref{eq:coneforL}. Since $a_1,a_2>0$ and $ma_1+na_2<0$, the remaining cases for the coefficients $m,n$ of $\al$ in \eqref{eq:condition} are $m<0,n\geq 0$ or $m\geq 0,n<0$. Without loss of generality, we assume $m<0,n\geq 0$. 
	
	Assume the hyperplane $\spn_\R\{\al\}$ is given by $P(A)$ for some $A\in E\bs\{0\}$. Then $\al\perp A$. Since $\{\al_1,\al_2\}$ is an $\R$-basis of $E$, we may choose its dual basis $\{\al^1,\al^2\}\subset E$ with respect to the inner product $(\cdot|\cdot)$ and assume 
	$$A=\la\al^1+\mu\al^2,\quad\al,\mu\in \R,\quad \mathrm{with}\quad  (\al|A)=\la m+\mu n=0.$$
	Since $m<0$ and $A\neq 0$, we have $\mu\neq 0$. Replace $A$ by $-A$, if necessary, we  assume $\mu>0$. Then 
	\begin{equation}\label{eq:intermediate1}
		\frac{\la}{\mu}+\frac{n}{m}=0.
	\end{equation}
	We claim that there exists two lattice point $\b,\b'\in L$ such that 
	\begin{enumerate}
		\item \label{(1)} $\b,\b'\in P^+(\ga)\cap L\subset M$.
		\item \label{(2)} $\{\al,\b\}$ is a basis of $L$, $\{\al,\b'\}$ is also a basis of $L$.
		\item  \label{(3)}$(A|\b)<0$ and $(A|\b')>0$ i.e., $\b,\b'$ are on different sides of the hyperplane $P(A)$. 
	\end{enumerate}
	See Figure~\ref{fig1} for an illustration of these vectors and hyperplanes in a standard rank-two lattice, where the shadow area represents $P^+(\ga)$.  
	
	If $\b$ and $\b'$ exists, then $\cone(\pm \al, \b)$ and $\cone(\pm \al,\b')$ are basis cones by \eqref{(2)}. By conditions \eqref{(1)}, \eqref{(3)}, and Lemma~\ref{lm:propertiesofcone}, $P^{\leq 0}(A)\cap L=\cone(\al,\b)\cup \cone (-\al,\b)\subset M$, and  $P^{\geq 0}(A)\cap L=\cone(\al,\b')\cup \cone (-\al,\b')\subset M$. Hence $L=M$. 
	
	\begin{figure}
		\centering
		\begin{tikzpicture}[scale=1.5]
			\clip (-4.2,-2.8) rectangle (4.2,2.8); 
			\draw[red, thick] (-4,2) -- (4,-2) node[below left] {$P(\ga)$};
			\draw[blue, dashed] (-4,-2) -- (4,2)node[above left] {$P(A)$};
			
			\foreach \x in {-4,-3.5,...,4} {
				\foreach \y in {-4,-3.5,...,4} {
					\fill[black] (\x,\y) circle (0.025);
				}
			}
			
			\fill[black] (0,0) circle (0.06);
			
			\draw[->, thick, blue] (0,0) -- (-0.5, 1) node[above right] {$A$};
			
			
			\draw[->, thick, blue] (0,0) -- (-1, -0.5) node[above left] {\small $\boldsymbol{\alpha}$};
			
			\draw[->, thick, red] (0,0) -- (0.5,1) node[right] {\small $\boldsymbol{\gamma}$};
			
			\draw[->, thick, orange] (0,0) -- (-1.5,1.5) node[above left] {\small $\boldsymbol{\beta}'$};
			
			\draw[->, thick, orange] (0,0) -- (2,0.5) node[above right] {\small $\boldsymbol{\beta}$};
			
			\fill [gray!30, fill opacity=0.4] plot    ({\x}, {-\x/2}) -- (4.5,-2.25) |- (-4,4);
			
		\end{tikzpicture}
		\caption{Illustration for the vectors and hyperplanes in Proposition~\ref{prop:monoidhyperplane} \label{fig1}}
	\end{figure}
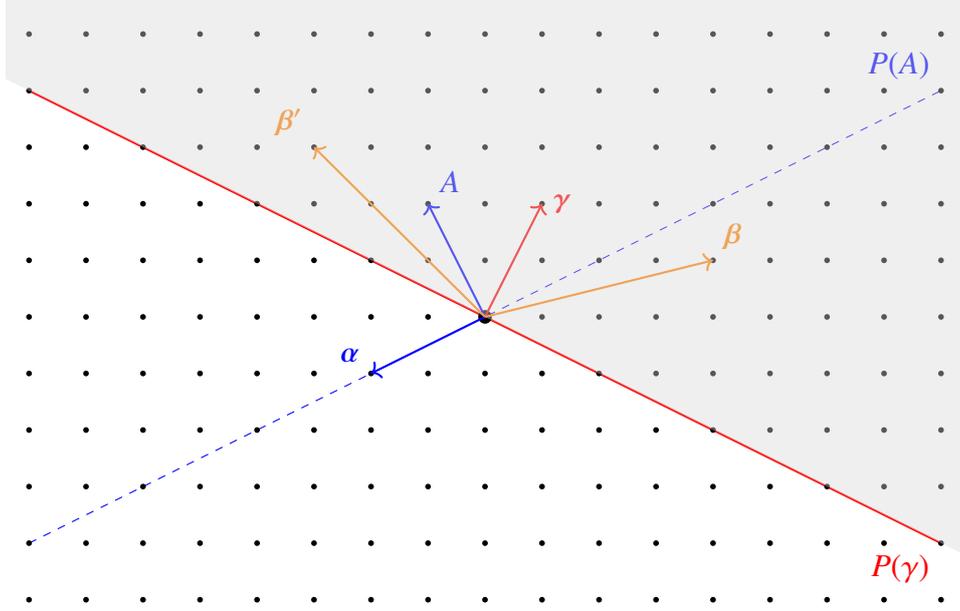

	To prove the existence of $\b$ and $\b'$, our idea is to first find some lattice points satisfying conditions \eqref{(1)} and \eqref{(2)}, then pick out the ones that satisfy condition \eqref{(3)}. We first prove the existence of $\b$. Let $x,y\in \Z$ be integral variables such that 
	\begin{equation}\label{eq:firsteq}
		\begin{bmatrix}
			\al\\\b
		\end{bmatrix}=\begin{bmatrix}
			m&n\\x&y
		\end{bmatrix}\begin{bmatrix}
			\al_1\\\al_2
		\end{bmatrix},\quad \mathrm{with}\quad my-nx=1,\quad \mathrm{and}\quad (\b|\ga)=xa_1+ya_2>0.
	\end{equation} 
	Any solution $\b=x\al_1+y\al_2$ to the system \eqref{eq:firsteq} would satisfy conditions \eqref{(1)} and \eqref{(2)} for $\b$. 
	
	Since $\al\in L$ is primitive, and we assumed $m<0$ and $n\geq 0$, then we have $\gcd(m,n)=1$ if $m,n\neq 0$; and $m=-1$ if $n=0$. In particular, there exists integers $(x_0,y_0)$ such that $my_0-nx_0=1$. Fix an integral solution $(x_0,y_0)$ to $my-nx=1$, the general integral solutions to the Diophantine equation $my-nx=1$ are given by 
	$$
	(x,y)=(x_0+tm,y_0+tn),\quad \mathrm{where}\quad  t\in \Z.
	$$
	Since $ma_1+na_2<0$ and $m<0$, we may choose a $t\in \Z_{<0}$ small enough so that  
	$$(x_0+tm)a_1+(y_0+tn)a_2=x_0a_1+y_0a_2+t(ma_1+na_2)>0,\quad \mathrm{and}\quad x_0+tm>0.$$
	Fix such a $t$, and let $(x_1,y_1)=(x_0+tm,y_0+tn)$. Then $\b=x_1\al_1+y_1\al_2$ is a solution to \eqref{eq:firsteq}, with $x_1>0$. Since $my_1-nx_1=1$, $x_1>0$ and $m<0$, we have   
	\begin{equation}\label{eq:intermediate2}
		\frac{y_1}{x_1}-\frac{n}{m}=\frac{1}{mx_1}. 
	\end{equation}
	Now we show $\b$ satisfies condition \eqref{(3)}. By \eqref{eq:intermediate1} and \eqref{eq:intermediate2}, together with the facts that $m<0$ and $x_1>0$, we have the following estimate:  
	$$\frac{(A|\b)}{\mu x_1}=\frac{(\la\al^1+\mu\al^2|x_1\al_1+y_1\al_2)}{\mu x_1}=\frac{\la}{\mu}+\frac{y_1}{x_1}=\frac{1}{mx_1}<0.$$
	Since $\mu x_1>0$, we have $(A|\b)<0$. 
	
	Now to show the existence of $\b'$, we consider another system similar to \eqref{eq:firsteq}: 
	\begin{equation}\label{eq:2ndeq}
		\begin{bmatrix}
			\al\\\b'
		\end{bmatrix}=\begin{bmatrix}
			m&n\\x&y
		\end{bmatrix}\begin{bmatrix}
			\al_1\\\al_2
		\end{bmatrix},\quad \mathrm{with}\quad my-nx=-1,\quad \mathrm{and}\quad (\b'|\ga)=xa_1+ya_2>0.
	\end{equation} 
	A similar argument shows that there exists integers $(x_2,y_2)$ such that 
	$$	\frac{y_2}{x_2}-\frac{n}{m}=-\frac{1}{mx_2} ,\quad x_2a_1+y_2a_2>0,\quad x_2>0.$$
	Then $\b'=x_2\al_1+y_2\al_2$ is a solution to \eqref{eq:2ndeq}, and 
	$$\frac{(A|\b')}{\mu x_2}=\frac{(\la\al^1+\mu\al^2|x_2\al_1+y_2\al_2)}{\mu x_2}=\frac{\la}{\mu}+\frac{y_2}{x_2}=-\frac{1}{mx_2}>0,$$
	since $m<0$ and $x_2>0$. Thus we have $(A|\b')>0$ since $\mu x_2>0$. 
\end{proof}

\begin{thm}\label{thm:classificationofparamonoid}
	Let $P$ be a parabolic-type submonoid of a rank-two even lattice $L$. Then it has two possible types: 
	\begin{enumerate}
		\item  $P=B$ is a rank-two Borel-type submonoid. In this case, $P=(\Z_{\geq 0}\al)\cup (P^+(\ga)\cap L)$, where $P(\ga)$ is the unique hyperplane such that $B\subset P^{\geq 0}(\ga)$ and $P(\ga)\cap L=\Z\al$ for some $\al\in L$ (possibly zero).
		\item  $P=P^{\geq 0}(\ga)\cap L$ for some hyperplane $P(\ga)\subset E$ such that $P(\ga)\cap L=\Z\al\neq \{0\}$.
	\end{enumerate}
	See Figures~\ref{fig2} and \ref{fig3} for an illustration in a standard rank-two lattice.
	We say that $P$ is of {\bf type-I} if $P=B$; and $P$ is of {\bf type-II} if $P=P^{\geq 0}(\ga)\cap L$ with $P(\ga)\cap L\neq \{0\}$. 
\end{thm}
\begin{figure}
	\centering
	\begin{tikzpicture}[scale=1.5]
		\clip (-4.2,-2.5) rectangle (4.2,2.5); 
		\draw[red, dashed] (0,0) -- (4,-2) node[below left] {$P(\ga)$};
		
		\foreach \x in {-4,-3,...,4} {
			\foreach \y in {-4,-3,...,4} {
				\fill[black] (\x,\y) circle (0.035);
			}
		}
		\foreach \x in {-3,-2,...,4} {
			\foreach \y in {2} {
				\fill[red] (\x,\y) circle (0.05);
			}
		}
		\foreach \x in {-1,0,...,4} {
			\foreach \y in {1} {
				\fill[red] (\x,\y) circle (0.05);
			}
		}
		\foreach \x in {0,1,...,4} {
			\foreach \y in {0} {
				\fill[red] (\x,\y) circle (0.05);
			}
		}
		
		\foreach \x in {2,3,4} {
			\foreach \y in {-1} {
				\fill[red] (\x,\y) circle (0.05);
			}
		}
		\fill[red](4,-2)  circle (0.05);
		
		\draw[->, thick, red] (0,0) -- (2, -1) node[above right] {$\al$};

		\draw[->, thick, red] (0,0) -- (0.5,1) node[right] {\small $\boldsymbol{\gamma}$};
		
		\fill [gray!30, fill opacity=0.4] plot    ({\x}, {-\x/2}) -- (4.5,-2.25) |- (-4,4);
		
		\node at (2,1.2) [above,red]{$P=B=(\Z_{\geq 0}\al)\cup (P^+(\ga)\cap L)$};
		
	\end{tikzpicture}
	\caption{Type-I parabolic-type submonoid} \label{fig2}
\end{figure}
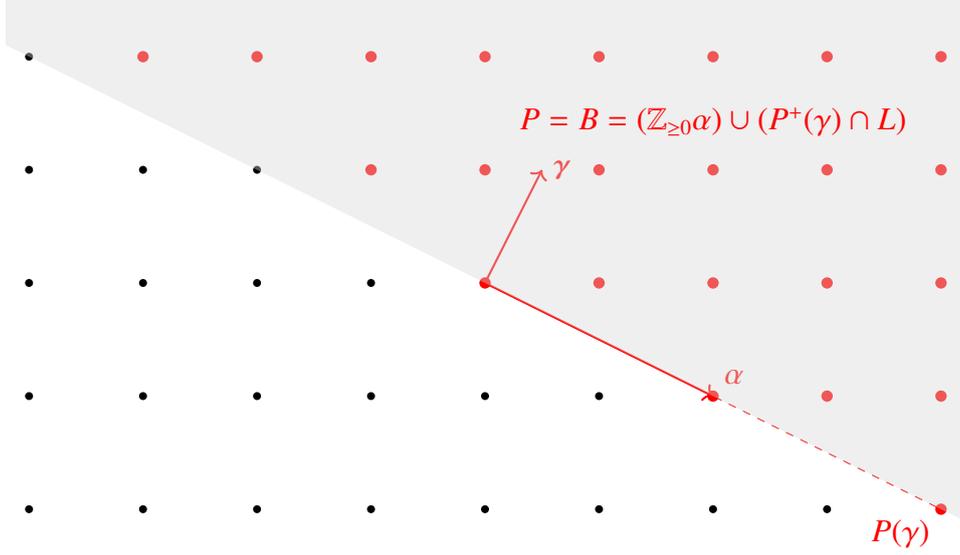

\begin{figure}
	\centering
	\begin{tikzpicture}[scale=1.5]
		\clip (-4.2,-2.5) rectangle (4.2,2.5); 
		\draw[red, dashed] (-4,2) -- (4,-2) node[below left] {$P(\ga)$};
		
		\foreach \x in {-4,-3,...,4} {
			\foreach \y in {-4,-3,...,4} {
				\fill[black] (\x,\y) circle (0.035);
			}
		}
		\foreach \x in {-4,-3,...,4} {
			\foreach \y in {2} {
				\fill[red] (\x,\y) circle (0.05);
			}
		}
		\foreach \x in {-2,-1,...,4} {
			\foreach \y in {1} {
				\fill[red] (\x,\y) circle (0.05);
			}
		}
		\foreach \x in {0,1,...,4} {
			\foreach \y in {0} {
				\fill[red] (\x,\y) circle (0.05);
			}
		}
		
		\foreach \x in {2,3,4} {
			\foreach \y in {-1} {
				\fill[red] (\x,\y) circle (0.05);
			}
		}
		\fill[red](4,-2)  circle (0.05);
		
		\draw[->, thick, red] (0,0) -- (2, -1) node[above right] {$\al$};
		
		\draw[->, thick, red] (0,0) -- (0.5,1) node[right] {\small $\boldsymbol{\gamma}$};
		
		\fill [gray!30, fill opacity=0.4] plot    ({\x}, {-\x/2}) -- (4.5,-2.25) |- (-4,4);
		
		\node at (2,1.2) [above,red]{$P=P^{\geq 0}(\ga)\cap L$};
		
	\end{tikzpicture}
	\caption{Type-II parabolic-type submonoid} \label{fig3}
\end{figure}
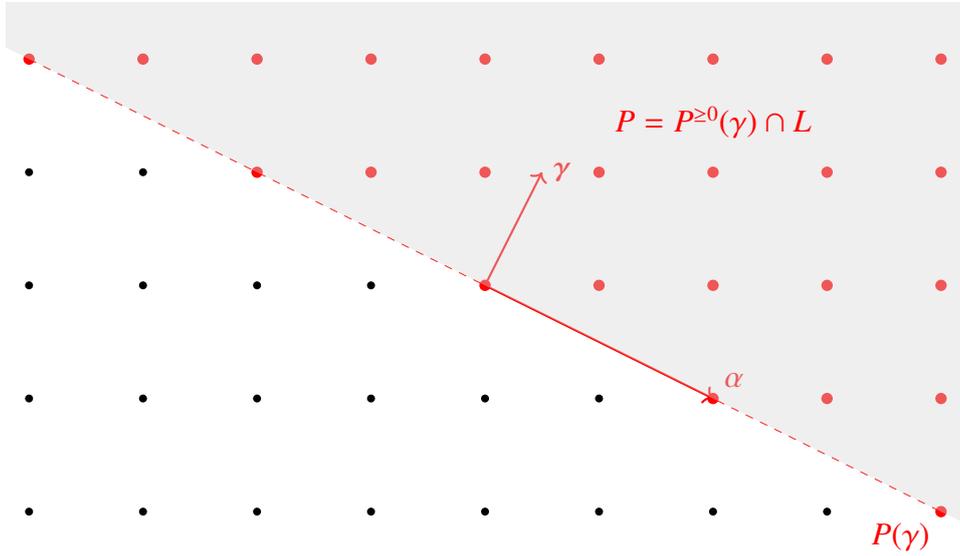
\begin{proof}
	
	By Definition~\ref{def:submonoidsofaPDElattice}, there exists a Borel-type submonoid $B$ such that $B\ssq P$. By Proposition~\ref{prop:basicBorel}, there exists a unique hyperplane $P(\ga)$ such that $B \subset P^{\geq 0}(\ga)$. Since $B\cap (-B)=\{0\}$ and $P^+(\ga)\cap P^{\leq 0}(\ga)=\emptyset$, we have $(-B)\ssq P^{\leq 0}(\ga)$ and $P^+(\ga)\cap L\subset B$. 
	
	Assume $P=B$ is of Type I. Note that $P(\ga)\cap L$ is a discrete additive subgroup in the one-dimensional Euclidean space $P(\ga)$. Then $P(\ga)\cap L=\Z\al$ for some $\al\in L$. Since $B\cup(-B)=L$, replace $\al$ by $-\al$ if necessary, we may assume $\al\in B$. Then
	the monoid $M=(\Z_{\geq 0}\al)\cup (P^+(\ga)\cap L)\ssq B$. Clearly, $M\cap(-M)=\{0\}=B\cap (-B)$, and $M\cup (-M)=L$. Hence $M=B$.

	Now assume $B\subsetneq P$. We claim that $P$ is of Type II. Indeed, by our argument above, $B=(\Z_{\geq 0}\al)\cup (P^+(\ga)\cap L)$. In particular $P^+(\ga)\cap L\subset P$. Let $\b\in P\bs B$. Then $\b$ is either in $\Z_{<0}\al$ or $P^{-}(\ga)\cap L$. If $\b\in P^{-}(\ga)\cap L$, write $\b=k\al'$ for some $k\geq 1$ and primitive vector $\al'\in P^{-}(\ga)\cap L$. Apply Proposition~\ref{prop:monoidhyperplane} to $M=P$, we have $P=L$, which contradicts the assumption that $P\subset L$ is a proper submonoid. Thus $\b\in \Z_{<0}\al$. Since $\al\in P$, it follows that $-\al\in P$ and so $(\Z\al)\cup (P^+(\ga)\cap L)=P^{\geq 0}(\ga)\cap L\ssq P$. Since $P$ contains no vectors in $P^-(\ga)\cap L$, we must have $P=P^{\geq 0}(\ga)\cap L$. 
\end{proof}


\section{Structural theory of $V_P$}\label{sec:3}
In this Section, we determine the structure of all parabolic-type VOAs $V_P$ associated to the rank-two parabolic-type submonoids $P$ based on our classification theorem for the parabolic-type monoids in the previous Section. 

\subsection{Basics of the lattice vertex operator algebras}\label{sec:basicsoflatticeVOA}
For the general definitions of vertex operator algebras (VOAs) and modules over VOAs, we refer to the classical texts \cite{FLM,FHL,DL,LL,FZ,Z}. Here we recall the notion of modules over a VOA. 

\begin{df}\label{df:modulesoverVOA}
	Let $V$ be a VOA. An {\bf admissible $V$-module} is a $\N$-graded vector space $M=\bigoplus_{n=0}^\infty M(n)$, equipped with a linear map $Y_M(\cdot,z):V\ra \End(M)[\![z,z^{-1}]\!],\ Y_M(a,z)=\sum_{n\in \Z} a_nz^{-n-1}$ called the {\bf module vertex operator}, satisfying 
	\begin{enumerate}
		\item (truncation property) For any $a\in V$ and $u\in M$, $Y_M(a,z)u\in M((z))$. 
		\item (vacuum property) $Y_M(\vac,z)=\Id_M$. 
		\item (Jacobi identity for $Y_M$) for any $a,b\in V$ and $u\in M$, 	
		\begin{equation}\label{eq:formalJacobi}
			\begin{aligned}
				z_0^{-1}\delta\left(\frac{z_1-z_2}{z_0}\right) &Y_M(a,z_1)Y_M(b,z_2)u-z_0^{-1}\delta\left(\frac{-z_2+z_1}{z_0}\right)Y_M(b,z_2)Y_M(a,z_1)u\\
				&=z_2^{-1}\delta\left(\frac{z_1-z_0}{z_2}\right) Y_M(Y(a,z_0)b,z_2)u.
			\end{aligned} 
		\end{equation}
		
		\item ($L(-1)$-derivative property) $Y_M(L(-1)a,z)=\frac{d}{dz} Y_M(a,z)$ for any $a\in V$.
		\item (grading property) For any $a\in V$, $m\in \Z$, and $n\in \N$,  $a_mM(n)\ssq M(n+\wt a-m-1)$. i.e., $\wt( a_m)=\wt a-m-1$. 
	\end{enumerate}
	We write $\deg v=n$ if $v\in M(n)$, and call it the {\bf degree} of $v$.  Submodules, quotient modules, and irreducible modules are defined in the usual categorical sense. Denote the category of admissible $V$-modules by $\mathsf{Adm}(V)$. $V$ is called {\bf rational} if the category $\mathsf{Adm}(V)$ is semisimple \cite{DLM1,Z}.
	
	An admissible $V$-module $M$ is called {\bf ordinary} if each degree-$n$ subspace $M(n)=M_{n+h}$ is a finite-dimensional eigenspace of $L(0)$ of eigenvalue $n+h$, where $h\in \C$ is called the {\bf conformal weight} of $M$. In particular, if we write $L(0)v=(\wt v)\cdot v$ for $v\in M(n)$, then $\wt v=\deg v+h$. 
	
	More generally, a {\bf weak} $V$-module is vector space $M$, together with a module vertex operator $Y_M(\cdot,z)$, satisfying conditions $(1),(2)$, and $(3)$ above. 
\end{df}

We briefly recall the construction of Heisenberg and lattice VOAs in \cite{FLM}. Some of the formulas here will be used later. 

Let $\h$ be a $r\geq 1$-dimensional vector space, equipped with a nondegenerate bilinear form $(\cdot|\cdot):\h\times \h\ra \C$. The affinization $ \hat{\h}=\h\otimes \C[t,t^{-1}]\op \C K$ has Lie bracket
\begin{equation}\label{2.3'}
	[h_1(m),h_2(n)]=m\delta_{m+n,0} (h_1|h_2) K,\quad [K,\hat{\h}]=0,\quad   h_1,h_2\in \h,\  m,n\in\Z,
\end{equation}
where $h\o t^m\in \hat{\h}$ is denoted by $h(m)$, together with a triangular decomposition $\hat{\h}=(\hat{\h})_{+}\op (\hat{\h})_{0}\op (\hat{\h})_{-}$, where $(\hat{\h})_{\pm}=\bigoplus_{n\in \Z_{\pm}} \h\o \C t^{n} $, and $(\hat{\h})_{0}=\h\o \C 1\op \C K$. Let $(\hat{\h})_{\geq 0}=(\hat{\h})_{+}\op (\hat{\h})_{0}$, which is a Lie sub-algebra of $\hat{\h}$. 

For each $\la\in \h$, let $e^\la$ be a formal symbol associated to $\la$. Then $\C e^\la$ is a module over $(\hat{\h})_{\geq 0}$, with $h(0)e^\la=(h|\la)e^\la$, $K.e^\la=e^\la$, and $h(n)e^\la=0$, for all $h\in \h$ and $n>0$. The induced module
\begin{equation}\label{eq:Heisenbergmod}
	M_{\hat{\h}}(1,\la)=U(\hat{\h})\o _{U(\hat{\h}_{\geq 0})} \C e^\la
\end{equation}
is an irreducible module over $\hat{\h}$. It has a basis 
\begin{equation}\label{2.4'}
	h_1(-n_1-1)\ds h_k(-n_k-1)e^\la: k\geq 0,\ h_1,\ds,h_k\in \h,\ n_{1}\geq \ds \geq n_k\geq 0 .
\end{equation}
It was proved in \cite{FLM,FZ} that $M_{\hat{\h}}(1,0)$ is a CFT-type simple VOA, with $\vac=e^0$ and $\om=\frac{1}{2}\sum_{i=1}^r u^i(-1)^2\vac$, where $\{u^1,\ds ,u^r\}$ is an orthonormal basis of $\h$, called the Heisenberg VOA of level-one, and $M_{\hat{\h}}(1,\la)$, with $\la\in \h$, are all the irreducible modules over the VOA $M_{\hat{\h}}(1,0)$ up to isomorphism. The (module) vertex operator is defined by normal ordering
\begin{equation}\label{eq:Heisenbergvertex}
	\begin{aligned}
		&Y( h_1(-n_1-1)\ds h_k(-n_k-1)\vac,z)={\tiny\begin{matrix}\circ \\\circ\end{matrix}}(\partial_{z}^{(n_1)}h_1(z))\dots (\partial_{z}^{(n_k)}h_k(z)){\tiny\begin{matrix}\circ \\\circ\end{matrix}}\\
		&\quad =\sum_{p,q,p+q=k} (\partial_{z}^{(n_{i_1})}h_{i_1}(z))_-\ds (\partial_{z}^{(n_{i_p})}h_{i_p}(z))_-(\partial_{z}^{(n_{j_1})}h_{j_1}(z))_+\ds (\partial_{z}^{(n_{j_q})}h_{j_q}(z))_+,
	\end{aligned}
\end{equation}
where the labels $(i_1,\ds,i_p,j_1,\ds,j_q)$ are permutations of $(1,2,\ds,k)$ and 
$$ \partial_z^{(n)}=\frac{1}{n!}\frac{\partial^n}{\partial z^n},\quad a(z)_-=\sum_{n<0} a_n z^{-n-1},\quad a(z)_+=\sum_{n\geq 0} a_n z^{-n-1},$$
for a field $a(z)=\sum_{n\in \Z} a_n z^{-n-1}$. 

Next, we recall the lattice VOAs. Let $L$ be an even lattice of rank $r\geq 1$ in a Euclidean space $E$, equipped with $\Z$-bilinear form $(\cdot|\cdot):L\times L\ra \Z$. Let $\h:=\C\otimes_{\Z} L$, extend $(\cdot|\cdot)$ to a $\C$-bilinear form $(\cdot|\cdot):\h\times \h\ra \C$. Let $\epsilon: L\times L\ra \<\pm 1\>$ be a $2$-cocycle of the abelian group $L$ such that $\epsilon(\al,\b)\epsilon(\b,\al)=(-1)^{(\al|\b)}$, for any $\al,\b\in L$. Write $\C^\ep[L]=\bigoplus_{\al\in L} \C e^{\al}$, where $e^{\al}$ is a formal symbol associated to $\al$ for each $\al\in L$ ($e^\al$ is denoted by $\iota(\al)$ in \cite{FLM}). Let 
\begin{equation}\label{eq:decofVL}
	V_L=M_{\hat{\h}}(1,0)\otimes_{\C} \C^\epsilon[L]=\bigoplus_{\al\in L} M_{\hat{\h}}(1,\al),
\end{equation}
where we identify $M_{\hat{\h}}(1,\al)$ with $M_{\hat{\h}}(1,0)\otimes \C e^\al$ for each $\al\in L$.
Define the lattice vertex operator $Y:V_L\ra \End(V_L)[[z,z^{-1}]]$ as follows
\begin{align}
	&Y(h(-1)\vac,z):=h(z)=\sum_{n\in \Z} h(n) z^{-n-1}\quad \left(h(n)e^\al:=0,\ \forall n>0,\ h(0)e^\la:=(h|\al) e^\al\right), \label{2.3}\\
	&Y(e^{\al},z):=E^{-}(-\al, z)E^{+}(-\al,z)e_{\al} z^{\al} \quad \left(z^\al (e^\b):=z^{(\al|\b)}e^\b,\ e_\al(e^\b):=\epsilon (\al,\b)e^{\al+\b}\right),\label{2.4}
\end{align}
where $\al,\b\in L$, $h\in \h$, and the operators $E^{\pm}$, $e^\al$, and $z^\al$ in \eqref{2.3}, \eqref{2.4} are given by
\begin{equation}\label{eq:latticevertex}
	E^{\pm}(-\al,z)=\exp\left(\sum_{n\in \Z_{\pm}}\frac{-\al(n)}{n}z^{-n}\right),\quad  e_\al(e^{\b})=\epsilon(\al,\b)e^{\al+\b},\quad z^\al(e^\b)=z^{(\al|\b)}e^{\b}.
\end{equation}
Moreover, the general lattice vertex operator is defined by normal ordering similar to \eqref{eq:Heisenbergvertex}
\begin{align*}
	&Y(h_1(-n_{1}-1)\dots h_k(-n_{k}-1)e^{\al},z):={\tiny\begin{matrix}\circ \\\circ\end{matrix}}(\partial_{z}^{(n_1)}h_1(z))\dots (\partial_{z}^{(n_k)}h_k(z))Y(e^{\al},z){\tiny\begin{matrix}\circ \\\circ\end{matrix}},\numberthis\label{2.5}\\
	&=\sum_{p,q,p+q=k} (\partial_{z}^{(n_{i_1})}h_{i_1}(z))_-\ds (\partial_{z}^{(n_{i_p})}h_{i_p}(z))_-E^{-}(-\al,z)e_\al z^\al E^{+}(-\al,z)(\partial_{z}^{(n_{j_1})}h_{j_1}(z))_+\ds (\partial_{z}^{(n_{j_q})}h_{j_q}(z))_+,
\end{align*}
where $k\geq 1$, $n_{1}\geq \ds\geq n_k\geq 0$, and $h_1,\dots h_k\in \h$.

It was proved in \cite[Appendix A.2]{FLM} that $V_L$ is a VOA, called the {\bf lattice VOA}. The Heisenberg VOA $M_{\hat{\h}}(1,0)$ is a subVOA of $V_L$ with the same Virasoro element $\om=\frac{1}{2}\sum_{i=1}^r u^i(-1)^2\vac$. Moreover, the lattice vertex operator $Y$ given by \eqref{2.3}--\eqref{2.5} are intertwining operators among the Heisenberg modules  $M_{\hat{\h}}(1,\al)$ in the decomposition \eqref{eq:decofVL} \cite{FLM,D,DL}. In particular, it satisfies 
\begin{equation}\label{eq:intertwin}
	Y(M_{\hat{\h}}(1,\al),z)M_{\hat{\h}}(1,\b)\subset M_{\hat{\h}}(1,\al+\b)[\![z,z^{-1}]\!],\quad \al,\b\in L.
\end{equation}

The irreducible modules over a lattice VOA were classified by Dong.

\begin{lm}\cite[Theorem 3.1]{D}\label{lm:VLmodules}
	Let $L$ be a even lattice in a Euclidean space $E$, and let $L^\circ=\{\la\in E: (\la|\al)\in \Z,\ \forall\al\in L \}$ be the dual lattice. Assume $L^\circ/L=\bigsqcup_{i=0}^r (L+\la_i)$ is the coset decomposition, then $\{V_{L+\la_i}:i=0,\ds, r\}$ are all the irreducible $V_L$-modules up to isomorphism. Moreover, $V_L$ is a rational VOA.
\end{lm}

The following fact follows from \eqref{eq:intertwin}.

\begin{lm}\cite[Proposition 3.2]{Liu25}.\label{prop2.3}
	Let $L$ be an even lattice in a Euclidean space $E$, $M\leq L$ be a submonoid, and $S\subset L$ be a sub-semigroup. Let $V_M=\bigoplus_{\al\in M} M_{\hat{\h}}(1,\al)$ and $V_S=\bigoplus_{\al\in S} M_{\hat{\h}}(1,\al)$. Then 
	\begin{enumerate}
		\item $(V_M,Y,\om ,\vac)$ is a CFT-type subVOA of the lattice VOA $(V_L,Y,\om ,\vac)$. 
		\item $(V_S,Y,L(-1))$ is a sub-vertex algebra without vacuum of $(V_L,Y,L(-1))$ \cite{HL}. If, furthermore, $S\subset M$ and $M+S\ssq S$, then $V_S$ is an ideal of $V_M$. 
	\end{enumerate}
\end{lm}

By choosing particular kinds of submonoid $M$ in a lattice $L$, we have different types of subVOAs in $V_L$. The following subVOAs in $V_L$ are our main objects to study in this paper. 

\begin{df}\label{def:typesofVOAs}
	Let $L=\Z\al_1\op \Z\al_2$ be a rank-two even lattice in a Euclidean space $E$. 
	\begin{enumerate}
		\item Let $P\leq L$ be a parabolic-type submonoid, see Definition~\ref{def:submonoidsofaPDElattice}. We call the subVOA $V_P=\bigoplus_{\al\in P} M_{\hat{\h}}(1,\al)$ of $V_L$ a {\bf rank-two parabolic-type VOA}. We say that the VOA $V_P$ is of {\bf type-I (resp. type-II)} if the rank-two submonoid $P$ is of  {\bf type-I (resp. type-II)}
		in Theorem~\ref{thm:classificationofparamonoid}. 
		\item Let $P(\ga)\subset E$ be a hyperplane. Then the submonoid $H=P(\ga)\cap L$ is a rank-one even lattice $\Z\al$ or $\{0\}$, see Theorem~\ref{thm:classificationofparamonoid}. When $H\neq \{0\}$, we call the subVOA 
		\begin{equation}\label{eq:VH}
			V_H:=\bigoplus_{n\in \Z} M_{\hat{\h}}(1,n\al)
		\end{equation}
		a {\bf rank-two hyperplane VOA}, see Figure~\ref{fig4} for an illustration. 
	\end{enumerate} 
\end{df}

\begin{remark}
	We observe the following facts from Definition~\ref{def:typesofVOAs}. 
	\begin{enumerate}
		\item The rank-two hyperplane VOA $V_H$ \eqref{eq:VH} is not isomorphic to the rank-one lattice VOA $V_{\Z\al}$ since $\dim \h=2$ and $M_{\widehat{\C\al}}(1,n\al)\neq M_{\hat{\h}}(1,n\al)$, see \eqref{2.4'} and \eqref{eq:decofVL}. 
		\item One can define parabolic-type VOAs in a higher rank lattice VOA, see \cite[Definition 3.4]{Liu25}. The notion of hyperplane VOA is also generalizable to the higher rank case. 
		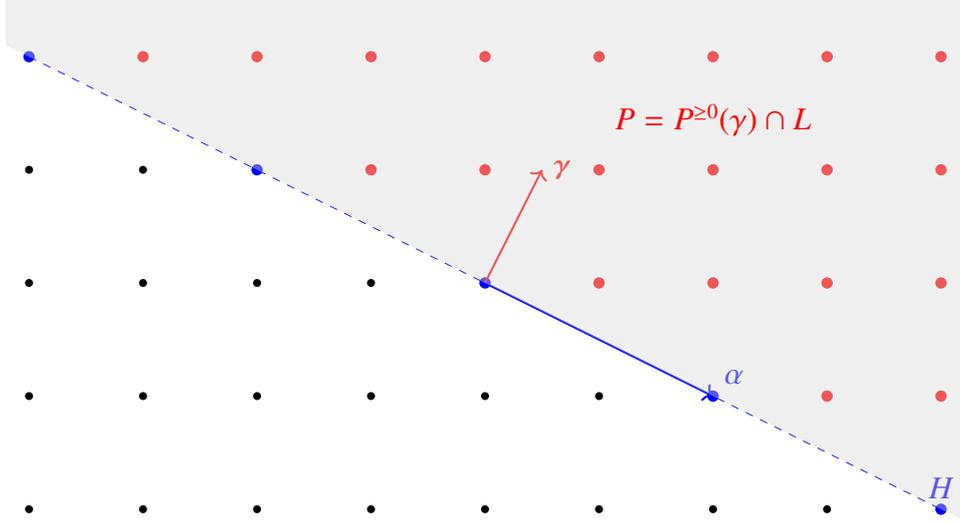
\begin{figure}
			\centering
			\begin{tikzpicture}[scale=1.5]
				\clip (-4.2,-2.5) rectangle (4.2,2.5); 
				\draw[blue, dashed] (-4,2) -- (4,-2) node[above ] {$H$};
				
				\foreach \x in {-4,-3,...,4} {
					\foreach \y in {-4,-3,...,4} {
						\fill[black] (\x,\y) circle (0.035);
					}
				}
				\foreach \x in {-3,-2,...,4} {
					\foreach \y in {2} {
						\fill[red] (\x,\y) circle (0.05);
					}
				}
				\foreach \x in {-1,0,...,4} {
					\foreach \y in {1} {
						\fill[red] (\x,\y) circle (0.05);
					}
				}
				\foreach \x in {1,2,...,4} {
					\foreach \y in {0} {
						\fill[red] (\x,\y) circle (0.05);
					}
				}
				
				\foreach \x in {3,4} {
					\foreach \y in {-1} {
						\fill[red] (\x,\y) circle (0.05);
					}
				}
				\fill[blue](4,-2)  circle (0.05);
				\fill[blue](2,-1)  circle (0.05);
				\fill[blue](0,0)  circle (0.05);
				\fill[blue](-2,1)  circle (0.05);
				\fill[blue](-4,2)  circle (0.05);
				
				\draw[->, thick, blue] (0,0) -- (2, -1) node[above right] {$\al$};
				
				\draw[->, thick, red] (0,0) -- (0.5,1) node[right] {\small $\boldsymbol{\gamma}$};
				
				\fill [gray!30, fill opacity=0.4] plot    ({\x}, {-\x/2}) -- (4.5,-2.25) |- (-4,4);
				
				\node at (2,1.2) [above,red]{$P=P^{\geq 0}(\ga)\cap L$};
				
			\end{tikzpicture}
			\caption{Hyperplane $H$ in type-II parabolic-type submonoid} \label{fig4}
		\end{figure}
	\end{enumerate}
\end{remark}

\subsection{Tensor product description of the rank-two hyperplane VOAs}

Recall the notion of tensor product VOAs introduced by Frenkel-Huang-Lepowsky. 
\begin{lm}\cite[Section 2.5, Proposition 4.7.2, Corollary 4.7.3]{FHL}\label{lm:tensor}
	Let $(V_1,Y_{V_1},\om_1,\vac_1)$ and $(V_2,Y_{V_2},\om_2,\vac_2)$ be two VOAs. The tensor product space $V_1\o V_2$ has a VOA structure with
	\begin{equation}\label{eq:tensor}
		Y_{V_1\o V_2}(a_1\o a_2,z)(b_1\o b_2)=Y_{V_1}(a_1,z)b_1\o Y_{V_2}(a_2,z)b_2,\quad a_1,b_1\in V_1, a_2,b_2\in V_2,
	\end{equation}
	$\vac_{V_1\o V_2}=\vac_1\o \vac_2$ and $\om_{V_1\o V_2}=\om_1\o \vac_2+\vac_1\o \om_2$.
	
	Furthermore, the irreducible ordinary modules over $V_1\o V_2$ are of the form $M_1\o M_2$, where $M_i$ are irreducible ordinary modules over $V_i$ for $i=1,2$, respectively. In particular, if $V_1$ and $V_2$ are both simple VOAs, then $V_1\o V_2$ is a simple VOA.  
\end{lm}

\begin{lm}\label{lm:3.5}
	Let $L$ be an even lattice in a Euclidean space $E$ of rank at least two. Let $\al,\b\in L\bs\{0\}$ such that $(\al|\b)=0$. Then for any $u=\al(-n_1-1)\ds \al(-n_r-1)e^{p\al}\in M_{\hat{\h}}(1,p\al)\subset V_L$, with $r\geq 0$, $n_1\geq \ds\geq n_r\geq 0$, and $p\in \Z$, the following relation holds in $\End(V_L)[\![z,z^{-1}]\!]$
	\begin{equation}\label{eq:commu}
		\b(m) Y(u,z)=Y(u,z)\b(m),\quad m\in \Z.
	\end{equation}
\end{lm}
\begin{proof}
	By \eqref{2.3'} we have $[\b(m),\al(n)]=0$ for any $m,n\in \Z$. Hence $\b(m)$ commutes with $(\partial^{(n)}\al(z))_\pm$ and $E^\pm(-p\al,z)$ for any  $n\geq 1$ and $p\in \Z$. Then by \eqref{2.5}, it suffices to show that $\b(m)$ commutes with the operator $e_{p\al}z^{p\al}$ in $\End(V_L)[\![z,z^{-1}]\!]$. 
	
	Indeed, if $m\neq 0$, by \eqref{2.3} and \eqref{2.4}, $\b(m)=\b(m)\o 1$ clearly commutes with $e_{p\al}z^{p\al}=1\o e_{p\al}z^{p\al}$ as operators on $M_{\hat{\h}}(1,\theta) \cong M_{\hat{\h}}(1,0)\otimes \C e^\theta\subset V_L$ for any $\theta\in L$. On the other hand, if $m=0$, then $\b(0)=1\o \b(0)$ as an operator on $M_{\hat{\h}}(1,0)\otimes \C e^\theta$. Since $(\b|p\al)=0$, we have
	\begin{align*}
		(\b(0)\circ (e_{p\al}z^{p\al}))e^\theta&=\b(0)\epsilon(p\al,\theta) z^{(p\al|\theta)} e^{p\al+\theta}=(\b|p\al+\theta)\epsilon(p\al,\theta) z^{(p\al|\theta)} e^{p\al+\theta}\\
		&=(\b|\theta)\epsilon(p\al,\theta) z^{(p\al|\theta)} e^{p\al+\theta}=( (e_{p\al}z^{p\al})\circ \b(0))e^\theta,
	\end{align*}
	for any $\theta\in L$. This proves \eqref{eq:commu}. 
\end{proof}


\begin{prop}\label{prop:VHtensor}
	The rank-two hyperplane VOA $V_H$ \eqref{eq:VH} is simple and is isomorphic to the tensor product of a rank-one Heisenberg VOA and a rank-one lattice VOA:
	$$M_{\widehat{\C \b}}(1,0)\o V_{\Z\al},$$
	where $\C\b=(\C\al)^\perp$ in $\h=\C\o_\Z L$ with respect to the extended $\C$-bilinear form $(\cdot|\cdot):\h\times \h\ra \C$. 
\end{prop}
\begin{proof}
	Since $\h=\C\al\op \C\b$ and $[\b(n),\al(m)]=0$ for all $m,n\in \Z$, the rank-two Heisenberg irreducible module $M_{\hat{\h}}(1,p\al)$ has a basis $\b(-n_1)\ds \b(-n_r)\al(-m_1)\ds \al(-m_s)e^{p\al}$, where $r,s\geq 0$, $n_1\geq \ds \geq n_r\geq 1$ and $m_1\geq \ds\geq m_s\geq 1$, see \eqref{2.4'}. Then we have an isomorphism of vector spaces:
	\begin{equation}\label{eq:defofvarphi}
		\begin{aligned}
			\varphi: M_{\hat{\h}}(1,p\al)&\ra M_{\widehat{\C\b}}(1,0)\o M_{\widehat{\C\al}}(1,p\al),\\
			\b(-n_1)\ds \b(-n_r)\al(-m_1)\ds \al(-m_s)e^{p\al}&\mapsto \b(-n_1)\ds \b(-n_r)\vac\o \al(-m_1)\ds \al(-m_s)e^{p\al},
		\end{aligned}
	\end{equation}
	which extends to a linear isomorphism 
	$$\varphi: V_H=\bigoplus_{p\in \Z} M_{\hat{\h}}(1,p\al)\xrightarrow{\sim} M_{\widehat{\C \b}}(1,0)\o \left(\bigoplus_{p\in \Z}M_{\widehat{\C\al}}(1,p\al) \right)=M_{\widehat{\C \b}}(1,0)\o V_{\Z\al}.$$
	Denote $M_{\widehat{\C \b}}(1,0)$ by $V_1$ and $V_{\Z\al}$ by $ V_2$ for short. We claim that
	\begin{equation}\label{eq:varphiiso}
		\varphi(Y_{V_H}(u,z)v)=Y_{V_1\o V_2}(\varphi(u),z)\varphi(v),
	\end{equation}
	for $u=\b(-n_1)\ds \b(-n_r)\al(-m_1)\ds \al(-m_s)e^{p\al}$ and $v=\b(-n'_1)\ds \b(-n'_{r'})\al(-m'_1)\ds \al(-m'_{s'})e^{q\al}$. 
	
	Indeed, by the commutativity of the lattice operators $\b(n)$ and $Y(\al(-m_1)\ds \al(-m_s)e^{p\al},z)$ as in Lemma~\ref{lm:3.5}, together with \eqref{2.5}, it is easy to see that 
	\begin{align*}
		Y_{V_H}(u,z)v&=\left({\tiny\begin{matrix}\circ \\\circ\end{matrix}}(\partial_{z}^{(n_1-1)}\b(z))\dots (\partial_{z}^{(n_r-1)}\b(z)){\tiny\begin{matrix}\circ \\\circ\end{matrix}}\right)\b(-n'_1)\ds \b(-n'_{r'})\\
		&\quad \cdot \left({\tiny\begin{matrix}\circ \\\circ\end{matrix}}(\partial_{z}^{(m_1-1)}\al(z))\dots (\partial_{z}^{(m_s-1)}\al(z))Y(e^{p\al},z){\tiny\begin{matrix}\circ \\\circ\end{matrix}}\right)\al(-m'_1)\ds \al(-m'_{s'})e^{q\al}.
	\end{align*}
	Then it follows from Lemma~\ref{lm:tensor} and \eqref{eq:defofvarphi}, together with the definition of vertex operators \eqref{eq:Heisenbergvertex} and \eqref{2.5} that
	\begin{align*}
		\varphi(Y_{V_H}(u,z)v)&=\left({\tiny\begin{matrix}\circ \\\circ\end{matrix}}(\partial_{z}^{(n_1-1)}\b(z))\dots (\partial_{z}^{(n_r-1)}\b(z)){\tiny\begin{matrix}\circ \\\circ\end{matrix}}\right)\b(-n'_1)\ds \b(-n'_{r'})\vac\\
		&\quad \o \left({\tiny\begin{matrix}\circ \\\circ\end{matrix}}(\partial_{z}^{(m_1-1)}\al(z))\dots (\partial_{z}^{(m_s-1)}\al(z))Y(e^{p\al},z){\tiny\begin{matrix}\circ \\\circ\end{matrix}}\right)\al(-m'_1)\ds \al(-m'_{s'})e^{q\al}\\
		&=Y_{V_1}(\b(-n_1)\ds \b(-n_r)\vac,z)\b(-n'_1)\ds \b(-n'_{r'})\vac \\
		&\quad \o Y_{V_2}(\al(-m_1)\ds \al(-m_s)e^{p\al},z)\al(-m'_1)\ds \al(-m'_{s'})e^{q\al}\\
		&=Y_{V_1\o V_2}(\b(-n_1)\ds \b(-n_r)\vac\o \al(-m_1)\ds \al(-m_s)e^{p\al},z)\\
		&\quad \cdot \b(-n'_1)\ds \b(-n'_{r'})\o \al(-m'_1)\ds \al(-m'_{s'})e^{q\al}\\
		&=Y_{V_1\o V_2}(\varphi(u),z)\varphi(v).
	\end{align*}
	This proves \eqref{eq:varphiiso}. Finally, let $\b_1=\b/\|\b\|$ and $\al_1=\al/\|\al\|$, where we define $\|h\|$ to be the principal branch of $(h|h)^{\frac{1}{2}}$ for $h\in \h$. Then $\{\b_1,\al_1\}$ is an orthonormal basis of $\h$ and $\om_{V_H}=\frac{1}{2}\b_1(-1)^2\vac+\frac{1}{2} \al_1(-1)^2\vac$. Moreover, $\frac{1}{2}\b_1(-1)^2\vac$ is the Virasoro element of $M_{\widehat{\C\b}}(1,0)$, and $\frac{1}{2}\al_1(-1)^2\vac$ is the Virasoro element of $V_{\Z\al}$. 
	Then by Lemma~\ref{lm:tensor} and \eqref{eq:defofvarphi}, we have
	\begin{align*}
		\varphi(\om_{V_H})=\frac{1}{2} \b_1(-1)^2\vac\o \vac+\vac\o\frac{1}{2} \al_1(-1)^2\vac=\om_{V_1}\o \vac+\vac\o \om_{V_2}=\om_{V_1\o V_2}. 
	\end{align*}
	Hence $\varphi: V_H\ra M_{\widehat{\C \b}}(1,0)\o V_{\Z\al}$ is an isomorphism of VOAs. Since the Heisenberg VOA $M_{\widehat{\C \b}}(1,0)$ and the lattice VOA $V_{\Z\al}$ are both simple, by Lemma~\ref{lm:tensor}, $V_H$ is a simple VOA. 
\end{proof}

\subsection{Structure of the rank-two parabolic-type VOAs $V_P$ and $C_1$-cofiniteness} 

Recall that the rank-two parabolic-type submonoids have the following classification, see Theorem~\ref{thm:classificationofparamonoid}.
\begin{equation}\label{eq:typesofP}
	\begin{aligned}
		&P=B=\Z_{\geq 0}\al \cup (P^+(\ga)\cap L),&& \mathrm{where}\quad P(\ga)\cap L=\Z\al,&&\mathrm{if}\ P\ \mathrm{is\ of\ type\ I},\\
		&P=P^{\geq 0}(\ga)\cap L,&& \mathrm{where}\quad P(\ga)\cap L=\Z\al\neq \{0\}, &&\mathrm{if}\ P\  \mathrm{is\ of\ type\ II}.
	\end{aligned}
\end{equation}

Now we determine the structures of rank-two parabolic-type VOAs.

\begin{thm}\label{thm:strucVP}
	Let $V_P$ be a rank-two parabolic-type VOA. Then there exists a unique maximum proper ideal $V^+\lhd V_P$ and a simple subVOA $V_H\leq V_P$ such that $V_P=V_H\op V^+$ as vector spaces and $V_P/V^+\cong V_H$ as VOAs. Moreover, we have the following characterization for the simple quotient VOA $V_H$: 
	\begin{enumerate}
		\item If $P=B$ is of type-I, then $V_H=M_{\hat{\h}}(1,0)$ is the rank-two Heisenberg VOA. 
		\item If $P$ is of type-II, then $V_H\cong M_{\widehat{\C \b}}(1,0)\o V_{\Z\al}$ is the hyperplane VOA (see \eqref{eq:VH} and Proposition~\ref{prop:VHtensor}), where $\b\in \h=\C\o_\Z L$ such that $\C \b=(\C\al)^\perp.$
	\end{enumerate}
	We call the subVOA $V_H$ the {\bf Cartan-part} of $V_P$. 
\end{thm}

\begin{proof}
	Assume $P$ is of type-I. Let $S:=P\bs\{0\}=\Z_{>0}\al \cup (P^+(\ga)\cap L)$ in view of  \eqref{eq:typesofP}. Clearly, $S\subset P$ is a sub-semigroup such that $P+S\ssq S$. Let $V^+:=V_S=\bigoplus_{\eta\in S}M_{\hat{\h}}(1,\eta)$. Then by Proposition~\ref{prop2.3},  $V^+\lhd V_P$ is an ideal. Moreover, we have $V_P=M_{\hat{\h}}(1,0)\op V^+$, with $V_H:=M_{\hat{\h}}(1,0)$ being a subVOA of $V_P$. Hence $V_P/V^+\cong V_H$ as VOAs. It remains that show that $V^+$ is the maximum proper ideal of $V_P$. 
	
	Indeed, let $J\lhd V_P$ be a proper ideal. It suffices to show $J\ssq V^+$. Since $J$ is closed under the action of the field $Y(h(-1)\vac,z)=h(z)=\sum_{n\in \Z} h(n)z^{-n-1}$, we have \begin{equation}\label{eq:J}
		J=\bigoplus_{\eta\in P}M_{\hat{\h}}(1,\eta)\cap J=\bigoplus_{\eta\in P}J_\eta.
	\end{equation}
	If $J_0=M_{\hat{\h}}(1,0)\cap J=0$, then $J\ssq V^+$, we are done. Otherwise, $J_0=M_{\hat{\h}}(1,0)$ since $M_{\hat{\h}}(1,0)$ is a simple VOA. Then $\vac\in M_{\hat{\h}}(1,0)\in J$ and $J=V_P$, which contradicts the properness of $J$. 
	
	Now assume $P$ is of type-II. Let $S:=P^+(\ga)\cap L$ and $V^+:=V_S$. Then by \eqref{eq:typesofP}, we have $P=S\sqcup H$, with $H=P(\ga)\cap L=\Z\al$ and $P+S\ssq S$. it follows that $V_P=V_H\op V^+$, $V_H\leq V_P$ is a subVOA, and $V^+\lhd V_P$ is an ideal. Hence $V_P/V^+\cong V_H$ as VOAs.  To show $V^+$ is the maximum proper ideal of $V_P$, we let $J\lhd V_P$ be a proper ideal. Then by \eqref{eq:J}, $J\cap V_H=\bigoplus_{\tau\in H} J_\tau$ is an ideal of the subVOA $V_H$. Since $\vac\in V_H$, we have $J\cap V_H\neq V_H$. But $V_H$ is a simple VOA by Proposition~\ref{prop:VHtensor}. Hence $J\cap V^+=0$ and $J= J\cap V^+\ssq V^+$. 
\end{proof}

\begin{remark}
	We note that the direct sum decomposition $V_P=V_H\op V^+$ in Theorem~\ref{thm:strucVP} can be viewed as a direct sum of a CFT-type VOA $V_H$ with a $V_H$-module $V^+$. However, this direct sum is {\bf not} the semi-direct product of the VOA $V_H$ with the module $V^+$ introduced by Li in \cite{L94} since $Y_{V_P}(V^+,z)V^+\neq 0$. 
\end{remark}

With the structure theorem for $V_P$, we discuss its $C_1$-cofiniteness. The subspace $C_1(V)\subset V$ was introduced by Li in \cite{L99}
\begin{equation}\label{eq:defofC1}
	C_1(V):=\spn\{a_{-1}b:a,b\in V_+ \}+\spn\{L(-1)v: v\in V_+\}.
\end{equation}
$V$ is called $C_1$-cofinite if $\dim V/C_1(V)<\infty$. 

A related notion is the strongly generation property of VOAs introduced by Kac in \cite{K}. A CFT-type VOA $V$ is called {\bf strongly generated} by a subset $U\ssq V$ if $V$ is spanned by 
\begin{equation}\label{2.13}
	u^{1}_{-n_1}u^2_{-n_2}\dots u^k_{-n_k}u,\quad \mathrm{where}\quad u^1,\dots ,u^k,u\in U,\quad n_{1}\geq n_{2}\geq \dots \geq n_{k}\geq 1.
\end{equation}
Karel and Li proved that a VOA $V$ is strongly generated by a finite-dimensional subspace $U\subset V$ if and only if $V$ is $C_1$-cofinite \cite{L99,KL}. 

\begin{prop}\label{prop:VHisC1}
	Let $V_P$ be a rank-two parabolic-type VOA. Then its Cartan-part subVOA $V_H$ is always $C_1$-cofinite. 
\end{prop}
\begin{proof}
	It is well-known that the Heisenberg VOA $M_{\hat{\h}}(1,0)$ and the lattice VOA $V_L$ are $C_1$-cofinite \cite{KL}. By Theorem~\ref{thm:strucVP}, to show $V_H$ is $C_1$-cofinite, it suffices to show if two CFT-type VOAs $V_1$ and $V_2$ are $C_1$-cofinite, then their tensor product $V_1\o V_2$ is also $C_1$-cofinite. By \eqref{eq:tensor}, it is clear that $C_1(V_1)\o V_2+V_1\o C_1(V_2)\ssq C_1(V_1\o V_2)$. Conversely, for $a^1\o a^2, b^1\o b^2\in (V_1\o V_2)_+$, since $L(-1)(a^1\o a^2)=L(-1)a^1\o a^2+a^1\o L(-1)a^2\in C_1(V_1)\o V_2+V_1\o C_1(V_2)$ and
	$$(a^1\o a^2)_{-1} (b^1\o b^2)=\sum_{i\in \Z} a^1_ib^1\o a^2_{-i-2}b^2\in C_1(V_1)\o V_2+V_1\o C_1(V_2).$$
	We have $C_1(V_1\o V_2)=C_1(V_1)\o V_2+V_1\o C_1(V_2)$ by \eqref{eq:defofC1}, and so 
	$$V_1\o V_2/C_1(V_1\o V_2)\cong (V_1/C_1(V_1))\o (V_2/C_1(V_2))$$
	is a finite-dimensional vector space. 
\end{proof}

Now we discuss the $C_1$-cofiniteness of the rank-two parabolic-type VOAs $V_P$. The following Lemma was proved in \cite[Lemma 5.1, Theorem 5.2, Remark 5.3]{Liu25}.
\begin{lm}\label{lm:ranktwocoincC1}
	Let  $L=\Z\al\op \Z\b$ be a rank-two even lattice, and $C=\mathrm{Cone}(\al,\b)=\Z_{\geq 0}\al+\Z_{\geq 0}\b$.
	\begin{enumerate}
		\item  If $(\al|\b)\geq 0$, then the conic-type VOA $V_C$ is $C_1$-cofinite. 
		\item If $(\al|\b)=-n$ for some $n\geq 1$, then the conic-type VOA $V_C$ is $C_1$-cofinite if the following condition holds:  Let $(\b|\b)=2k$ and $(\al|\al)=2\ell$, with $k\geq \ell\geq 1$, the integers $n,k,\ell$ need to satisfy 
		\begin{equation}\label{eq:conditionforranktwoC_1}
			n^2+\ell^2-4\ell k\leq 0.
		\end{equation}
		In particular, if $(\al|\al)=2$, then $V_C$ is $C_1$-cofinite. 
	\end{enumerate}
\end{lm}

\begin{lm}\label{lm:sumC1}
	Let $V$ be a CFT-type VOA, and $V_1,V_2$ be two CFT-type subVOAs that have the same vacuum and Virasoro elements with $V$ itself. If $V=V_1+V_2$ and $V_1$, $V_2$ are both $C_1$-cofinite, then $V$ is also $C_1$-cofinite. 
\end{lm}
\begin{proof}
	By \eqref{eq:defofC1}, it is clear that $C_1(V_1)+C_2(V_2)\ssq C_1(V)$. Then there exists canonical surjective linear maps
	$$V_1/C_1(V_1)\op V_2/C_1(V_2)\twoheadrightarrow V/(C_1(V_1)+C_1(V_2))\twoheadrightarrow V/C_1(V).$$
	Hence $\dim (V/C_1(V))\leq \dim (V_1/C_1(V_1))+\dim (V_2/C_1(V_2))<\infty$. 
\end{proof}

\begin{thm}\label{thm:C1forVP}
	Let $V_P$ be a rank-two parabolic-type VOA. 
	\begin{enumerate}
		\item Assume $P=B$ is of type-I, see Theorem~\ref{thm:classificationofparamonoid}. Then $V_P$ is {\bf not} $C_1$-cofinite.
		\item Assume $P=P^{\geq 0}(\ga)\cap L$ is of type-II, with $P(\ga)\cap L=\Z\al\neq \{0\}$, see Theorem~\ref{thm:classificationofparamonoid}.
		
		If there exists $\b\in P$ such that $\{\al,\b\}$ is a basis of $L$, and condition \eqref{eq:conditionforranktwoC_1} is satisfied for the ordered basis $\{\al,\b\}$ or $\{\b,\al\}$, then $V_P$ is $C_1$-cofinite. In particular, if $(\al|\al)=2$ or $(\b|\b)=2$, then $V_P$ is $C_1$-cofinite. 
	\end{enumerate}
\end{thm}
\begin{proof}
	The proof of (1) is similar to the proof of \cite[Proposition 5.7]{Liu25} for the non-$C_1$-cofiniteness of the standard Borel-type subVOA $V_B$ of $V_{A_2}$. We briefly sketch the idea. Assume  $P=B=(\Z_{\geq 0}\al)\cup (P^+(\ga)\cap L)$, with $P(\ga)\cap L=\Z\al$. Choose $\b\in P$ such that $\{\al,\b\}$ is a basis of $L$. Then $\b-m\al\in P$ for any $m\geq 1$. One can show that $\{e^{\b-m\al}\in V_P:m=1,2,\ds\}$ is not contained in $C_1(V_P)$, and the image of this set is linearly independent in $V_P/C_1(V_P)$. We omit the further details. 
	
	Now assume $P=P^{\geq 0}(\ga)\cap L$ is of type-II, and there exists $\b\in P$ such that $\{\al,\b\}$ is a basis of $L$ and condition \eqref{eq:conditionforranktwoC_1} is satisfied for $\{\al,\b\}$. It follows from \eqref{eq:coneforL} that $P=\cone(\al,\b)\cup \cone(-\al,\b)$.  Let $C_1=\cone(\al,\b)$ and $C_2=\cone(-\al,\b)$. Then $V_P=V_{C_1}+V_{C_2}$ as a VOA, see Lemma~\ref{prop2.3}. Then $V_{C_1}$ is $C_1$-cofinite by Lemma~\ref{lm:ranktwocoincC1} (2). On the other hand, since $(\al|-\b)=n\geq 0$, then $V_{C_2}$ is $C_1$-cofinite by Lemma~\ref{lm:ranktwocoincC1} (1). Hence $V_P=V_{C_1}+V_{C_2}$ is $C_1$-cofinite by Lemma~\ref{lm:sumC1}. 
\end{proof}

\section{Representation theory of $V_P$}\label{sec:4}

In this Section, we present a detailed investigation of the representation theory of the rank-two parabolic-type VOA $V_P$. In particular, we classify all the irreducible modules over $V_P$ and determine the fusion rules among them. Our main results indicate that the representation theory of the parabolic-type VOA $V_P$ is governed by its Cartan-part $V_H$, generalizing the representation theoretical results of the rank-one Borel-type VOAs \cite[Theorem 6.13]{Liu25}. 

\subsection{Nil ideal in the Zhu's algebra $A(V_P)$}
We want to show that the maximum ideal $V^+$ of $V_P$ (see Theorem~\ref{thm:strucVP}) acts as $0$ on any irreducible $V_P$-module. So that the irreducible modules of $V_P$ are in one-to-one correspondence with the irreducible modules over its Cartan-part $V_H$. 
However, this is not obvious from the definition of $V_P$. We invoke the important tool of Zhu's algebra $\A=A(V)$ \cite{Z} to prove this claim. 

First, we recall the definition and basic properties $A(V)$. Let $(V,Y,\vac,\om)$ be a VOA. For homogeneous elements $a,b\in V$, define
\begin{align}
	a\circ b&:=\Res_{z=0} Y(a,z)b\frac{(1+z)^{\wt a}}{z^2}=\sum_{j\geq 0} \binom{\wt a}{j} a_{j-2}b\label{3.1},\\
	a\ast b&:=\Res_{z=0}Y(a,z)b\frac{(1+z)^{\wt a}}{z}=\sum_{j\geq 0} \binom{\wt a}{j} a_{j-1}b\label{3.2}.
\end{align}
Let $O(V)=\spn\{a\circ b:a,b\in V\}$, and let $A(V):=V/O(V)$. By \cite[Theorem 2.1.1]{Z}, $O(V)$ is a two-sided ideal with respect to $\ast$, and $A(V)$ is an associative algebra with respect to $\ast$, with the unit element $[\vac]$. The following formulas can be found in \cite[Lemma 2.1.3]{Z}: 
\begin{align}
	&	a\ast b\equiv \Res_{z=0} Y(b,z)a\frac{(1+z)^{\wt b-1}}{z} \pmod{O(V)},\label{3.3}\\
	&	\Res_{z=0} Y(a,z)b\frac{(1+z)^{\wt a+n}}{z^{2+m}}\equiv 0\pmod{O(V)},\label{3.5}
\end{align}
where $a,b\in V$ are homogeneous, and $m\geq n\geq 0$. Recall the following fundamental facts about Zhu's algebra.

\begin{lm}\cite[Theorem 2.2.2]{Z}\label{lm:1-1}
	Let $W$ be a weak $V$-module. Then the subspace
	$$\Om(W)=\spn\{w\in W: a_kw=0,\ \deg(a_k)<0\}$$
	is a left $A(V)$-module via $o:A(V)\ra \End(\Om(W)), [a]\mapsto o(a)=a_{\wt a-1}$. 
	In particular, the degree-zero subspace $M(0)$  of an admissible $V$-module $M=\bigoplus_{n=0}^\infty M(n)$ a left $A(V)$-module.
	
	Moreover, $M\mapsto M(0)$ is a one-to-one correspondence between the isomorphism class of irreducible admissible $V$-modules  and the isomorphism class of irreducible $A(V)$-modules. In particular, $A(V)$ is a finite-dimensional semisimple algebra if $V$ is rational.
\end{lm}

\begin{lm}\cite[Proposition 1.4.2]{FZ}\label{lm:ideal}
	Let $V$ be a VOA and $I\lhd V$ be an ideal. Then $A(I)=(I+V)/O(V)$ is a two-sided ideal of $A(V)$. 
	
	In particular, if $V_P=V_H\op V^+$ is a rank-two parabolic-type VOA, then $A(V^+)=\spn\{[a]:a\in V^+\}$ is a two-sided ideal of $A(V_P)$. 
\end{lm}

\begin{prop}\label{prop:nilideal}
	Let $V_P$ be a rank-two parabolic-type VOA. Then $A(V^+)$ is a nil ideal of the Zhu's algebra $A(V_P)$. Furthermore, if $V_P$ is $C_1$-cofinite, see Theorem~\ref{thm:C1forVP} for a sufficient condition, then $A(V^+)$ is a nilpotent ideal of $A(V_P)$. 
\end{prop}
\begin{proof}
	By Theorem~\ref{thm:strucVP}, $A(V^+)=\spn\{[M_{\hat{\h}}(1,\b)]:\b\in S \}$, where $S=P\bs\{0\}$ if $P$ is of type-I, and $S=P^+(\ga)\cap L$ if $P$ is of type-II. 
	
	We claim that $[M_{\hat{\h}}(1,2\b)]=0$ in $A(V_P)$ for any $\b\in S$.
	Indeed, since $0\notin S$, we may assume $(\b|\b)=2N$ for some $N\geq 1$. Then $\wt (e^\b)=N$. By \eqref{2.4} and \eqref{eq:latticevertex}, 
	\begin{align*}
		e^\b_{-2N-1}e^\b&=\Res_{z=0}z^{-2N-1}(E^-(-\b,z)E^+(-\b,z)e_\b z^{\b})e^\b \\
		&= \Res_{z=0} z^{-2N-1} \exp\left(\sum_{n<0} \frac{-\b(n)}{n}z^{-n}\right)\epsilon(\b,\b) z^{2N} e^{2\b}=e^{2\b}.
	\end{align*}
	Similarly, we can show that $e^\b_{-m}e^\b=0$ for $m\leq 2N$.
	Then it follows from \eqref{3.5} that 
	\begin{align*}
		0\equiv \Res_{z=0} Y(e^\b,z)e^\b \frac{(1+z)^{\wt (e^\b)}}{z^{2N+1}}=\sum_{j\geq 0}\binom{N}{j} e^\b_{-2N-1+j} e^\b=e^{2\b}\pmod{O(V_P)}. 
	\end{align*}
	Moreover, it follows from \eqref{3.1} and \eqref{3.3} that $h(-n-2)u\equiv -h(-n-1)u\pmod{O(V_P)}$ and $[h(-1)u]=[u]\ast [h(-1)\vac]$ in $A(V_P)$, for any $u\in V_P$, $h\in \h$, and $n\geq 0$. Then we have
	\begin{align*}
		[h^1(-n_1-1)\ds h^r(-n_r-1)e^{2\b}]=(-1)^{n_1+\ds +n_r}[e^{2\b}]\ast [h^1(-1)\vac]\ast\ds \ast [h^r(-1)\vac]=[0],
	\end{align*} 
	for any $h^i\in \h$ and $n_1\geq \ds \geq n_r\geq 0$, since $[e^{2\b}]=0$. This shows $[M_{\hat{\h}}(1,2\b)]=0$ in $A(V_P)$ for any $\b\in S$.
	
	Now let $\al\in [M_{\hat{\h}}(1,\b)]\subset A(V^+)$. By \eqref{eq:intertwin} and \eqref{3.2} we have 
	\begin{align*}
		\al\ast \al \in [M_{\hat{\h}}(1,\b)]\ast [M_{\hat{\h}}(1,\b)]\subset [M_{\hat{\h}}(1,2\b)]=[0]. 
	\end{align*}
	Hence each element in $A(V^+)$ is nilpotent, and so $A(V^+)$ is a nil ideal of $A(V_P)$.
	
	If $V_P$ is $C_1$-cofinite, then it follows from \cite[Theorem 3.1]{Liu21} that $A(V_P)$ is a Noetherian algebra. Hence the nil ideal $A(V^+)$ is nilpotent by Levitzky's theorem \cite{Le45}. 
\end{proof}

\begin{coro}\label{coro:nilaction}
	Let $V_P$ be a rank-two parabolic-type VOA, and let $M=\bigoplus_{n=0}^\infty M(n)$ be an irreducible admissible $V_P$-module. Then the action of the ideal $A(V^+)$ on $M(0)$ via $o:A(V_P)\ra \End(M(0))$ is zero. 
\end{coro}
\begin{proof}
	By Lemma~\ref{lm:1-1}, $M(0)$ is an irreducible module over $A(V_P)$. Since $A(V^+)\lhd A(V_P)$ is a nil ideal by Proposition~\ref{prop:nilideal}, it is contained in the Jacobson radical $J(A(V_P))$ which acts as zero on any irreducible $A(V_P)$-module. 
\end{proof}

\subsection{Classification of irreducible modules over $V_P$}

By Theorem~\ref{thm:strucVP}, $V_P=V_H\op V^+$, with $V_H\leq V_P$ and $V^+\lhd V_P$. Let $(M,Y_{M})$ be an irreducible admissible $V_H$-module. Define 
\begin{equation}\label{eq:tY}
	\tilde{Y}_{M}:V_P\ra \End(M)[\![z,z^{-1}]\!],\quad \tY_M|_{V_H}:=Y_M,\ \tY_M|_{V^+}:=0. 
\end{equation}
It is clear that $(M,\tY_M)$ is an irreducible admissible $V_P$-module. 

\begin{lm}\label{lm:irrvanishing}
	Let $V_P$ be a rank-two parabolic-type VOA and $M=\bigoplus_{n=0}^\infty M(n)$ be an irreducible admissible $V_P$-module. Then the action of the ideal $V^+$ on $M$ must be zero. 
\end{lm}
\begin{proof}
	Since $V^+$ is an ideal of $V_P$, its action on $M$ gives rise to an admissible submodule $N=V^+.M$ which has spanning elements 
	\begin{equation}\label{eq:4.5}
		a^1_{n_1}\ds a^r_{n_r}w,\quad \mathrm{where}\quad r\geq 1,\ a^i\in V^+,\ n_i\in \Z,\ w\in M.
	\end{equation}
	We need to show that $N=0$. 
	
	Indeed, suppose $N\neq 0$. Since $M$ is irreducible, then $N=M$. We can also fix a nonzero element $u\in M(0)$ and write
	$M=\spn\{b_mu:b\in V_P,m\in \Z\}$, see \cite[Proposition 4.5.7]{LL}. Then it follows from \eqref{eq:4.5} that $M(0)=N(0)$ has spanning elements
	\begin{equation}\label{eq:4.6}
		a^1_{n_1}\ds a^r_{n_r}b_mu,\quad \mathrm{where}\quad r\geq 1,\ a^i\in V^+,\ b\in V_P,\ \sum_{i=1}^r \deg (a^i_{n_i})+\deg (b_m)=0.
	\end{equation}
	We use induction on $r\geq 1$ to show that elements of the form \eqref{eq:4.6} are $0$. 
	
	Consider the base case $r=1$. Note that $b_mu=0$ if $\deg (b_m)<0$. Hence we may assume $\deg (a^1_{n_1})\leq 0$ in $a^1_{n_1}b_mu$. By Corollary \ref{coro:nilaction} we have $o(V^+)M(0)=0$.
	If $\deg (a^1_{n_1})=0$, then $\deg (b_m)=0$ and 
	$a^1_{n_1}b_mu=o(a^1)o(b)u=0$ since $a^1\in V^+$ and $o(b)u\in M(0)$. If $\deg (a^1_{n_1})<0$, then 
	$a^1_{n_1}u=0$ and 
	$$a^1_{n_1}b_mu=[a^1_{n_1},b_m]u=\sum_{j\geq 0}\binom{n_1}{j} (a^1_jb)_{n_1+m-j}u=\sum_{j\geq 0}\binom{n_1}{j} o(a^1_jb)u=0,$$
	since  $\deg (a^1_{n_1})+\deg (b_m)=0$ and $a^1_jb\in V^+$ for all $j\geq 0$. Now assume \eqref{eq:4.6} holds for smaller $r$. If $\deg (a^1_{n_1})=\ds =\deg (a^r_{n_r})=0$ then $a^1_{n_1}\ds a^r_{n_r}b_mu=o(a^1)\ds o(a^r)o(b)u=0$. Otherwise, there exists $i$ such that $\deg (a^i_{n_i})<0$. Then 
	\begin{align*}
		a^1_{n_1}\ds a^r_{n_r}b_mu&=\sum_{k=i}^r\sum_{j\geq 0}\binom{n_i}{j} a^1_{n_1}\ds \widehat{a^i_{n_i}}\ds (a^i_ja^k)_{n_i+n_k-j}\ds a^r_{n_r}b_mu \\
		&+\sum_{j\geq 0}\binom{n_i}{j} a^1_{n_1}\ds \widehat{a^i_{n_i}}\ds  a^r_{n_r}(a^i_jb)_{n_i+m-j}u.
	\end{align*}
	Since each summand on the right-hand-side has shorter length, $a^1_{n_1}\ds a^r_{n_r}b_mu=0$ by the induction hypothesis. Thus, $M(0)=N(0)=0$, which is a contradiction. 
\end{proof}

We also need the following result about representation of tensor product associative algebras. 

\begin{lm}\label{lm:tensormodule}
	Let $A$ be a finite-dimensional semisimple algebra over $\C$, with irreducible modules $S^1,\ds,S^r$ up to isomorphism, and let $B=\C[x]\o A$ be the tensor product associative algebra. For any $\la\in \C$, let $e^\la$ be a formal symbol. Then $\C e^\la\o S^i$ is an irreducible $B$-module with
	\begin{equation}\label{eq:tensoraction}
		\rho: B\ra \End(\C e^\la\o S^i),\quad \rho(f(x)\o a)(e^\la\o u)=f(\la)\cdot e^\la\o a.u,\quad a\in A,u\in S^i.
	\end{equation}
	Furthermore, any irreducible $B$-module is isomorphic to $\C e^\la\o S^i$ for some $\la\in \C$ and $1\leq i\leq r$. 
\end{lm}
\begin{proof}
	Clearly, $\C e^\la\o S^i$ is an irreducible $B$-module via \eqref{eq:tensoraction}. 
	However, it is not straightforward that any irreducible $B$-module has the form $\C e^\la\o S^i$, since the subalgebra $\C[x]$ of $B$ is not finite-dimensional. 
	
	Since $A=M_{n_1}(\C)\times \ds \times M_{n_r}(\C)$ by the Artin-Wedderburn's theorem, it suffices to show any irreducible module over the tensor product algebra $\C[x]\o M_n(\C)\cong M_n(\C[x])$ is isomorphic to $\C e^\la \o S$, where $S=\C^n$ is the standard irreducible $M_n(\C)$-module.

	Let $M\neq 0$ be an irreducible $M_n(\C[x])$-module. Since $M$ is a module over the semisimple subalgebra $M_n(\C)$, there exists a nonzero copy of $S=\C^n\subset M$. Let $N\ssq M$ be the $\C[x]$-submodule generated by $S$. i.e., $N=\spn\{f(x).u\in M:f(x)\in \C[x],\ u\in S \}$. Clearly, $N\leq M$ is a nonzero $M_n(\C[x])$-submodule, and so $N=M$. Moreover, $\C[x]\o S\ra M$, $f(x)\o u\mapsto f(x).u$ is an epimorphism of $M_n(\C[x])$-modules. Then $M=(\C[x]\o S)/J$, where $J$ is a maximal proper submodule of $\C[x]\o S\cong \C[x]^n=\C[x]\op\ds \op \C[x]$. We claim there exists $\la\in \C$ such that
	\begin{equation}\label{eq:maxJ}
		J=\spn\left\{ \begin{bmatrix}
			f_1(x)\\\vdots \\f_n(x)
		\end{bmatrix}\in \C[x]^n: f_i(x)\in \<x-\la\>,\ \forall 1\leq i\leq n \right\}.
	\end{equation}
	Indeed, if there exists two elements $w_1=[\ds f_i(x)\ds ]^t$ and $w_2=[\ds g_j(x)\ds]^t$ in $J$ such that $\gcd(f_i(x),g_j(x))=1$, assume $a(x)f_i(x)+b(x)g_j(x)=1$ for some $a(x),b(x)\in \C[x]$, then $(a(x)E_{1i}).w_1+(b(x)E_{1j}).w_2=[1,0,\ds,0]^t\in J$, where $E_{ij}\in M_n(\C)$ is the matrix unit. It follows that $\C[x]^n=M_n(\C[x]).[1,0,\ds,0]^t\ssq J$, which contradicts the properness of $J$. Hence there exists a nonzero polynomial $p(x)$ such that for any $[f_1(x),\ds, f_n(x)]^t\in J$, we have $p(x)|f_i(x)$ for all $i$. It follows that $J\ssq \<p(x)\>^n=\<p(x)\>\op \ds \op \<p(x)\>$. By the maximality of $J$, we have $p(x)=x-\la$ for some $\la\in \C$, and $J=\<p(x)\>^n$. This proves \eqref{eq:maxJ}. 
	
	Now it follows from \eqref{eq:maxJ} that $M\cong \C[x]^n/\<x-\la\>^n\cong (\C[x]/\<x-\la\>)^n\cong \C e^\la\o \C^n$ as a module over $M_n(\C[x])\cong \C[x]\o M_n(\C)$. 
\end{proof}


\begin{thm}\label{thm:classificationirr}
	Let $V_P=V_H\op V^+$ be a rank-two parabolic-type VOA and $M=\bigoplus_{n=0}^\infty M(n)$ be an irreducible admissible $V_P$-module. Then $M$ is an irreducible $V_H$-module on which $V^+$ acts as zero. In particular, any irreducible admissible $V_P$-module is ordinary. Furthermore,
	\begin{enumerate}
		\item If $V_P$ is of type-I, then 
		\begin{equation}\label{eq:typeIirr}
			\{(M_{\hat{\h}}(1,\la),\tY_M):\la\in \h\}
		\end{equation} are all the irreducible $V_P$-modules up to isomorphism, where $\tY_M$ is defined by \eqref{eq:tY};
		\item  If $V_P$ is of type-II, with $P(\ga)\cap L=\Z\al$, $(\al|\al)=2N$, and  $\C \b=(\C\al)^\perp$, then 
		\begin{equation}\label{eq:typeIIirr}
			\{(L^{(\mu,i)}=M_{\widehat{\C\b}}(1,\mu)\o V_{\Z\al+\frac{i}{2N}\al},\tY_M): \mu\in \C\b,\ i=0,1,\ds, 2N-1 \}
		\end{equation}
		are all the irreducible $V_P$-modules up to isomorphism, where $\tY_M$ is defined by \eqref{eq:tY}.
	\end{enumerate}
\end{thm}
\begin{proof}
	By Lemma~\ref{lm:irrvanishing}, $M$ is an irreducible admissible $V_H$-module on which $V_+$ acts as $0$. The non-trivial part is to show that $M$ is ordinary. i.e., $\dim M(n)<\infty$ for all $n\geq 0$.
	
	If $V_P$ is of type-I, we have $V_H=M_{\hat{\h}}(1,0)$ by  Theorem \ref{thm:strucVP}. Then $M(0)$ is an irreducible module over $A(M_{\hat{\h}}(1,0))\cong \C[x,y]$ by Lemma ~\ref{lm:1-1} and \cite[Theorem 3.1.2]{FZ}. By Hilbert's Nullstellensatz, $M(0)\cong \C e^\la$ for some $\la\in \h$. Since the irreducible Heisenberg module $M_{\hat{\h}}(1,\la)$ is a Verma $\hat{\h}$-module by its construction \eqref{eq:Heisenbergmod}, we have $M_{\hat{\h}}(1,\la)\cong M$, and so $M$ is ordinary. This also proves \eqref{eq:typeIirr}. 
	
	If $V_P$ is of type-II, then $V_H\cong M_{\widehat{\C\b}}(1,0)\o V_{\Z\al}$ by Theorem \ref{thm:strucVP}. It was proved in  \cite[Lemma 2.8]{DMZ94} that $A(V_H)\cong A(M_{\widehat{\C\b}}(1,0))\o A(V_{\Z\al})\cong \C[x]\o A(V_{\Z\al})$ as associative algebras. Then by Lemma~\ref{lm:VLmodules} and Lemma~\ref{lm:1-1}, $A(V_{\Z\al})$ is a semisimple algebra, and $U= M(0)$ is an irreducible module over $A(V_H)\cong \C[x]\o A(V_{\Z\al})$. 
	Hence $\dim U< \infty$ and $U\cong \C e^\mu \o S$ by Lemma~\ref{lm:tensormodule}, where $\mu \in \C\b$ and $S$ is an irreducible $A(V_{\Z\al})$-module. Since $(\Z\al)^\circ /\Z\al=\bigsqcup_{i=0}^{2N-1}(\Z\al+\frac{i}{2N}\al)$, then $S$ is the bottom degree of  $V_{\Z\al+\frac{i}{2N}\al}$ for some $0\leq i\leq 2N-1$ by Lemma~\ref{lm:VLmodules}. 
	In particular, $U=M(0)$ is isomorphic to the bottom degree of the irreducible ordinary $V_H$-module $W=M_{\widehat{\C\b}}(1,\mu)\o V_{\Z\al+\frac{i}{2N}\al}$. Note that both $M$ and $W$ are the unique irreducible quotients of the generalized Verma module $\bar{M}(U)$ associated to the irreducible $A(V_H)$-module $U$ \cite[Theorems 6.2, 6.3]{DLM1}. Hence $M\cong W=M_{\widehat{\C\b}}(1,\mu)\o V_{\Z\al+\frac{i}{2N}\al}$ is an ordinary module. 
	This also proves \eqref{eq:typeIIirr}. 
\end{proof}

\begin{remark}
	In the argument of the type-II case of Theorem \ref{thm:classificationirr}, we cannot directly apply Frenkel-Huang-Lepowsky's result for irreducible ordinary modules over tensor product VOA (see Lemma~\ref{lm:tensor}) to the irreducible admissible $V_H$-module $M=\bigoplus_{n=0}^\infty M(n)$ without knowing $M$ is ordinary. In fact, there exists VOAs (e.g. affine VOAs of admissible level) whose irreducible admissible modules are not necessarily ordinary. Hence we have to show $M$ is an ordinary $V_H$-module first. In fact, \eqref{eq:typeIIirr} agrees with  Lemma~\ref{lm:tensor} once we have the ordinarity of $M$. 
\end{remark}


\section{Strongly unital property of the Cartan-part subVOA $V_H$}\label{sec:5}

In this Section, we prove that the Cartan-part subVOA $V_H$ of a type-II parabolic-type VOA $V_P$ \eqref{eq:VH} satisfies the strongly unital property \cite{DGK23,DGK24}. These are new natural examples of CFT-type $C_1$-cofinite irrational simple VOA that satisfies the strongly unital property adding to the Heisenberg VOA example.

\subsection{The mode transition algebras $\fA_d$ associated to a VOA}
We refer to \cite{FZ,FBZ04,DGK23} for the definition of universal enveloping algebra  $\scrU=\scrU(V)$ of a VOA $V$. 

$\scrU=\bigoplus_{d\in \Z}\scrU_d$ is a canonically seminorm graded associative algebra. It has left and right neighborhoods at $0$:
\begin{equation}\label{eq:nebors}
	\rN^n_\L \scrU=\scrU\cdot \scrU_{\leq -n}\quad \mathrm{and}\quad  	\rN^n_\sR \scrU=\scrU_{\geq n}\cdot \scrU,
\end{equation}
such that $	\rN^{n+1}_\L \scrU_0=\sum_{j \geq n+1} \scrU_j\cdot \scrU_{-j}=	\rN^{n+1}_\sR \scrU_0$ for any $n\geq 0$, see \cite[Appendix A]{DGK23}.
The Zhu algebra $\A=A(V)$ has the following identification as an associative algebra. 
\begin{equation}
	\A\cong\frac{\scrU_0}{\rN^{1}_\L \scrU_0}=\frac{\scrU_0}{\sum_{j \geq 1} \scrU_j\cdot \scrU_{-j}}
\end{equation}
Denote the category of $\A=A(V)$-modules by $\mathsf{Mod}(\A)$. 

On the other hand, let $W$ be a left module over the associative algebra $\scrU$  equipped with discrete topology. If the module action map $\scrU\times W\ra W$ is continuous, then $W$ must be a filtered module
\begin{equation}
	W=\bigcup_{n=0}^\infty\Om_n(W),\quad \mathrm{where}\quad  \Om_n(W)=\{w\in W: (\rN^{n+1}_\L\scrU).w=0\}. 
\end{equation}

We call such a $\scrU$-module $W$ {\bf exhaustive} or continuously discrete.
In particular, any admissible $V$-module $M=\bigoplus_{n=0}^\infty M(n)$ is naturally an exhaustive $\scrU$-module, with $M(n)\ssq \Om_n(M)$ for all $n$, see Definition~\ref{df:modulesoverVOA}. Denote the category of exhaustive left $\scrU$-modules by $\mathsf{ExMod}(\scrU)$. Then $\mathsf{Adm}(V)$ is a full subcategory of  $\mathsf{ExMod}(\scrU)$.

It was observed in \cite{DGK23} that $\Om_n(W)$ agrees with the space of degree-$n$ highest-weight vectors in the weak $V$-module $W$ introduced in \cite{DLM1}. i.e., 
$$\Om_n(W)=\{w\in W: a_{k}w=0,\ \deg(a_{k})\leq -n-1\}.$$
In particular, 
\begin{equation}\label{eq:Om}
	\Om: \mathsf{Adm}(V)\ra \mathsf{Mod}(\A),\quad W\mapsto \Om_0(W)
\end{equation}
is a functor between abelian categories. 

On the other hand, Dong-Li-Mason's generalized Verma module functor $\bar{M}(\cdot)$ \cite[Theorem 6.2]{DLM1} can be identified with the following (left) induced module functor:
\begin{equation}\label{eq:PhiL}
	\Phi^\L:  \mathsf{Mod}(\A)\ra  \mathsf{Adm}(V)\subset \mathsf{ExMod}(\scrU),\quad S\mapsto \Phi^\L(S)=(\scrU/	\rN^1_\L \scrU)\o _{\scrU_0} S. 
\end{equation}
The functors $\Om$ and $\Phi^\L$ in \eqref{eq:Om} and \eqref{eq:PhiL} form an adjoint pair $\Phi^\L\dashv \Om$ between abelian categories
\begin{equation}\label{eq:adjointequi}
	\Phi^\L:\mathsf{Mod}(\A)\rightleftarrows\mathsf{Adm}(V) :\Om,
\end{equation}
see \cite[Proposition 3.1.2]{DGK23}. The following Lemma is useful for our later discussion.
\begin{lm}\label{lm:PhiOmfortensor}
	Let $V=V_1\o V_2$ be a tensor product VOA. Then the functors $\Phi^\L$ and $\Om$, with $\A=A(V_1\o V_2)$ and $\scrU=\scrU(V_1\o V_2)$, have the following properties
	\begin{enumerate}
		\item $\Phi^\L(U)=\Phi^\L_{V_1}(\Phi^\L_{V_2}(U))=\Phi^\L_{V_2}(\Phi^\L_{V_1}(U))$ for any $ U\in \mathsf{Mod}(\A)$.
		\item $\Om(W)=\Om_{V_1}(W)\cap \Om_{V_2}(W)$ for any $W\in \mathsf{Adm}(V).$
	\end{enumerate}
\end{lm}
\begin{proof}
	Both (1) and (2) are straightforward consequences of the fact that $(a_1\o a_2)_{[n]}$ is a linear combination of operators of the form $(a_1\o \vac)_{[p]}(\vac\o a_2)_{[q]}=( \vac\o a_2)_{[q]}(a_1\o \vac)_{[p]}$, see \cite[eq. (4.7.3)]{FHL}. More precisely, given any weak $V=V_1\o V_2$-module $W$, by \eqref{eq:tensor}, 
	\begin{equation}\label{eq:moduletensor}
		(a_1\o a_2)_{[n]}.w=\sum_{i\in \Z} \left((a_1)_{[i]}\o (a_2)_{[n-i-1]} \right).w=\sum_{i\in \Z} (a_1\o \vac)_{[i]}.((\vac\o a_2)_{[n-i-1]}.w). 
	\end{equation}
	Clearly, \eqref{eq:moduletensor} makes $W=\Phi^\L_{V_1}(\Phi^\L_{V_2}(U))$ an admissible $V$-module generated by $U$. Then there exists a canonical epimorphism $\pi: \Phi^\L(U)\ra \Phi^\L_{V_1}(\Phi^\L_{V_2}(U))$. On the other hand, since $V_1,V_2$ are subVOAs of $V_1\o V_2$ and the actions of $V_1$ and $V_2$ are commutative on a $V$-module, there exists a $V_2$-module homomorphism $\iota: \Phi^\L_{V_2}(U)\ra \Phi^\L(U)$ such that $\Im \iota\ssq \Om_{V_1}(\Phi^\L(U))$. Note that $\iota$ also extends to a $V_1$-module homomorphism $\tau: \Phi_{V_1}^\L(\Phi_{V_2}^\L(U))\ra \Phi^\L(U)$ by \cite[Theorem 6.2]{DLM1}. By \eqref{eq:moduletensor}, $\tau$ is also a $V$-homomorphism and is clearly an inverse of $\pi$. This proves (1).  
	
	Clearly, $\Om(W)\ssq \Om_{V_1}(W)\cap \Om_{V_2}(W)$ since $V_1$ and $V_2$ are subVOAs of $V$. On the other hand, assume $\deg((a_1\o a_2)_{[n]})=\wt a_1+\wt a_2-n-1<0$, then for any $w\in \Om_{V_1}(W)\cap \Om_{V_2}(W)$ and any $i\in \Z$, we have 
	$$(a_1\o \vac)_{[i]}.((\vac\o a_2)_{[n-i-1]}.w)=(\vac\o a_2)_{[n-i-1]}.((a_1\o \vac)_{[i]}.w)=0 $$
	since either $\deg ((a_1\o \vac)_{[i]})=\wt a_1-i-1<0$ or $\deg ((\vac\o a_2)_{[n-i-1]})=\wt a_2-n+i<0$. Hence $w\in \Om(W)$ by \eqref{eq:moduletensor}. 
\end{proof}

Using the right neighborhoods in \eqref{eq:nebors}, one can define a right induced module functor 
\begin{equation}\label{eq:PhiR}
	\Phi^\sR: \mathsf{Mod}^r(\A)\ra  \mathsf{ExMod}^r(\scrU),\quad U\mapsto \Phi^\sR(U)=U\o_{\scrU_0} (\scrU/\rN^1_\sR \scrU).
\end{equation}
Note that the neighborhood $\rN^1_\L \scrU$ (resp. $\rN^1_\sR \scrU$) contains all the negatively (resp. positively) graded subspaces of $\scrU$, see  \eqref{eq:nebors}. Hence the quotient algebras in \eqref{eq:PhiL} and \eqref{eq:PhiR} have the gradations 
$$\scrU/	\rN^1_\L \scrU=\bigoplus_{n=0}^\infty(\scrU/	\rN^1_\L \scrU)_n,\quad \scrU/\rN^1_\sR \scrU=\bigoplus_{m=0}^\infty (\scrU/\rN^1_\sR \scrU)_{-m}.$$
\begin{df}\cite{DGK23,DGK24}\label{df:MTA}
	The {\bf mode transition algebra} $\fA=\fA(V)$ associated to a VOA $V$ is given by the following vector space
	\begin{equation}
		\fA=\Phi^\sR(\Phi^\L(\A))=(\scrU/	\rN^1_\L \scrU)\o _{\scrU_0} \A \o_{\scrU_0} (\scrU/\rN^1_\sR \scrU).
	\end{equation}
	Write $\fA_{n,-m}=(\scrU/	\rN^1_\L \scrU)_{n}\o _{\scrU_0} \A\o _{\scrU_0} (\scrU/\rN^1_\sR \scrU)_{-m}$, then $\fA=\bigoplus_{n,m\in \N} \fA_{n,-m}$. The product on $\fA$ is denoted by $\star$, which satisfies 
	$$\fA_{i,-j}\star \fA_{k,-l}\ssq \delta_{j,k}\fA_{i,-l},$$
	see \cite[Appendix B, Definition 3.2.1]{DGK23} for more details of the product. $\fA_d=\fA_{d,-d}$ is an associative algebra under $\star$, called the {\bf $d$-th mode transition algebra} of $V$. 
	
	An element $1_d\in \fA_d$ is called a {\bf strong unit} if it satisfies 
	$$1_d\star \mathfrak{a}=\mathfrak{a}\quad \mathrm{and}\quad \mathfrak{b}\star 1_d=\mathfrak{b},\quad \mathfrak{a}\in  \fA_{d,0},\ \mathfrak{b}\in  \fA_{0,-d},$$
	where $\fA_d\cong \fA_{d,0}\o_\A \fA_{0,-d}$ and 
	\begin{equation}\label{eq:fAd0}
		\fA_{d,0}=(\scrU/	\rN^1_\L \scrU)_{d}\o _{\scrU_0} \A\quad \mathrm{and}\quad  \fA_{0,-d}=\A \o_{\scrU_0} (\scrU/\rN^1_\sR \scrU)_{-d}
	\end{equation}
	are right and left $\A$-modules, respectively. 
	$\fA_d$ is said to be {\bf strongly unital} if it admits a strong identity element. We say that the VOA $V$ satisfies the {\bf strongly unital condition for mode transition algebras}, or simply $V$ is strongly unital, if its mode transition algebras $\fA_d$ are strongly unital for all $d\in \N$. 
\end{df}

\begin{lm}\cite[Theorem 4.0.10]{DGK24}\cite[Theorem 6.1.1]{DGK23}\label{lm:egofstronglyunital}
	Any rational VOA satisfies the  strongly unital condition for mode transition algebras. Moreover, any Heisenberg VOA $V=M_{\hat{\h}}(k,0)$ of rank $r\geq 1$ satisfies the  strongly unital condition for mode transition algebras.
\end{lm}

\subsection{Strongly unital property of  $V_H$} 
The following theorem that characterizes the strongly unital condition was proved in an ongoing work \cite{GL26}. 

\begin{lm}\label{lm:equicondforstronglyunital}
	Let $V$ be a $C_1$-cofinite VOA such that $\fA_{d,0}$ is a projective right $\A$-module for all $d\geq 0$. The following conditions are equivalent:
	\begin{enumerate}
		\item The VOA $V$ satisfies the  strongly unital condition for mode transition algebras. 
		\item The adjoint pair $\Phi^\L\dashv \Om$ in \eqref{eq:adjointequi} is an adjoint equivalence.
	\end{enumerate}
	In particular, any admissible $V$-module $W$ is a generalized Verma module $\Phi^\L(\Om(W))$. 
\end{lm}

As a consequence, it is clear that the rank-two parabolic-type VOA $V_P$ does not satisfy the strongly unital property since it is not a simple VOA and so it is not a generalized Verma module.

\begin{thm}\label{thm:stronglyuntialforVP}
	Let $V_P=V_H\op V^+$ be a rank-two parabolic-type VOA. Then its Cartan-part subVOA $V_H$ satisfies the  strongly unital condition for mode transition algebras. 
\end{thm}
\begin{proof}
	By Proposition~\ref{prop:VHisC1}, $V_H$ is $C_1$-cofinite. If $V_P$ is of type-I, then $V_H\cong M_{\hat{\h}}(1,0)$ by Theorem~\ref{thm:strucVP}, and the conclusion follows from Lemma~\ref{lm:egofstronglyunital}. 
	
	Assume $V_P$ is of type-II, with $V_H\cong M_{\widehat{\C\b}}(1,0)\o V_{\Z\al}=V_1\o V_2$ as VOAs, see Theorem~\ref{thm:strucVP}. Both $V_1=M_{\hat{\h}}(1,0)$ and $V_2=V_{\Z\al}$ are strongly unital by Lemma~\ref{lm:egofstronglyunital} and Lemma~\ref{lm:VLmodules}. Moreover, it follows from the definition of universal enveloping algebra and tensor product vertex operators \eqref{eq:tensor} that 
	\begin{align*}
		&	\scrU(V_H)_0\cong \scrU(V_1)_0\o_\C \scrU(V_2)_0,\\
		&\scrU(V_H)_d\o _{\scrU(V_H)_0}W\cong \bigoplus_{i,j\geq 0,i+j=d}(\scrU(V_1)_i\o_\C \scrU(V_2)_j)\o_{\scrU(V_H)_0}W,
	\end{align*}
	for any left $A(V_H)$-module $W$. Then by \eqref{eq:fAd0}, 
	$$\fA_{d,0}(V_H)\cong \bigoplus_{i,j\geq 0,\ i+j=d}\fA_{i,0}(V_1)\o_\C \fA_{j,0}(V_2)$$
	as a right $A(V_H)\cong A(V_1)\o_\C A(V_2)$-module. Since $\fA_{i,0}(V_1)$ (resp. $\fA_{j,0}(V_2)$) is a projective right $A(V_1)$ (resp. $A(V_2)$) module, respectively, it is clear that $\fA_{d,0}(V_H)$ is a projective right $A(V_H)$-module. Now let $U\in \mathsf{Mod}(A(V_H))$. By Lemma~\ref{lm:PhiOmfortensor}, 
	\begin{align*}
		\Om\Phi^\L(U)&=\Om(\Phi^\L_{V_1}(\Phi^\L_{V_2}(U)))= \Om_{V_1}\Phi^\L_{V_1}(\Phi^\L_{V_2}(U))\cap \Om_{V_2}\Phi^\L_{V_2}(\Phi^\L_{V_1}(U))\\
		&=\Phi^\L_{V_2}(U)\cap \Phi^\L_{V_1}(U)
	\end{align*}
	since the unit maps  for $V_1$ and $V_2$ are both bijective. Suppose there exists some element 
	$$w=\sum_i \b(-n_1)\ds \b(-n_k)u^i=\sum_{j} a^1_{m_1}\ds a^r_{m_r}v^j\in \left(\Phi^\L_{V_2}(U)\cap \Phi^\L_{V_1}(U)\right)\bs U, $$
	where $n_1\geq \ds \geq n_k\geq 1$, $a^1,\ds, a^r\in V_{\Z\al}$, and $u^i,v^j\in U$. Choose the index $i$ such that $n_1$ is largest among all the summands $\b(-n_1)\ds \b(-n_k)u^i$. Then there exists a nonzero scalar $c\in \C$ such that 
	$$c\cdot u^i=\b(n_k)\ds \b(n_1).w=\sum_{j} a^1_{m_1}\ds a^r_{m_r}\b(n_k)\ds \b(n_1).v^j=0,$$
	in view of \eqref{2.3'}. Then $u^i=0$, which is a contradiction. Thus, $\Om\Phi^\L(U)=\Phi^\L_{V_2}(U)\cap \Phi^\L_{V_1}(U)=U$.
	On the other hand, given $W\in \mathsf{Ext}(\scrU(V_H))$, since the action of $V_1$ and $V_2$ are commutative on $W$, $\Om_{V_2}(W)$ is a $V_1$-submodule in $W$, and $\Om(W)=\Om_{V_1}(\Om_{V_2}(W))$. Since $V_2$ is strongly unital,  $\Phi^\L_{V_2}\Om_{V_2}(W)= W$ by Lemma~\ref{lm:equicondforstronglyunital}. Then $\Om_{V_1}(\Phi^\L_{V_2}\Om_{V_2}(W))= \Om_{V_1}(W)$ and 
	$$\Phi^\L(\Om(W))=\Phi^\L_{V_1}\Om_{V_1}(\Phi^\L_{V_2}\Om_{V_2}(W))= \Phi_{V_1}^\L\Om_{V_1}(W)= W,$$
	in view of Lemma~\ref{lm:PhiOmfortensor}. Hence $V_H=V_1\o V_2$ satisfies the  strongly unital condition for mode transition algebras by Lemma~\ref{lm:equicondforstronglyunital}.
\end{proof}


\end{document}